\documentclass[12pt,reqno]{amsart}

\usepackage{amsmath,amsfonts,amssymb,amscd,amsthm,amsbsy,epsf}

 \newtheorem{thm}{Theorem}[section]

 \theoremstyle{definition}
 
 \theoremstyle{remark}
 
 \numberwithin{equation}{section}

 \newcommand{\N}{\mathbb{N}}
 \newcommand{\Z}{\mathbb{Z}}
 
 \newcommand{\C}{\mathbb{C}}

\begin{document}

\title[Identities Hurwitz and L-functions]
{Identities for the Hurwitz zeta function, Gamma function, and $L$-functions}

\author[Michael O. Rubinstein]
{Michael O.\ Rubinstein\\ \\
Pure Mathematics \\University of Waterloo\\200 University Ave W\\Waterloo, ON, N2L 3G1\\Canada}




\begin{abstract}
We derive several identities for the Hurwitz and Riemann zeta functions, the
Gamma function, and Dirichlet $L$-functions. They involve a sequence of polynomials
$\alpha_k(s)$ whose study was initiated in~\cite{R}.
The expansions given here are practical and can be used for the high precision
evaluation of these functions, and for deriving formulas for special
values. We also present a summation formula and use it to generalize a formula of Hasse.
\end{abstract}

\maketitle

\section{Introduction}

Let $s \in \C$ and define $\alpha_k(s)$ to be the Taylor coefficients given by:
\begin{equation}
    \label{eq:1}
    \left(\frac{-\log (1-t)}{t}\right)^{s-1} = \sum_0^\infty \alpha_k (s) t^k\ ,\qquad |t|<1.
\end{equation}
Then $\alpha_0(s)=1$, $\alpha_1(s)=(s-1)/2$,
and, in \cite{R}, the following formula is proved:
\begin{equation}
    \label{eq:recur}
    \alpha_{k+1}(s)
    = \frac1{k(k+1)(k+2)} \sum_{j=1}^k
    \frac{\alpha_j(s)j(k+k^2+s(2k+2-j))}{(k-j+1)(k-j+2)} \ , \qquad k\geq 1.
\end{equation}
The above recursion shows that, for $k \geq 1$, $\alpha_k(s)/(s-1)$ is a polynomial in $s$ with
positive rational coefficients, and allows one to obtain the bound~\cite{R}:
\begin{equation}
    \label{eq:alpha bound}
    |\alpha_k(s)| \leq c_s \frac{(1+\log(k+1))^{|s|+1}}{k+1}\ ,
\end{equation}
where
\begin{equation}
    c_s = \frac{|s-1|}{|s|+1} (|s|+2) 2^{|s|+1}.
\end{equation}
This bound was used in~\cite{R} in the derivation of the following formulas for the Riemann
zeta and Gamma functions, valid for all $s \in \C$ not a pole of the relevant function:
\begin{equation}
    \label{eq:gamma2 intro}
    \Gamma (s) = \sum_0^\infty \frac{\alpha_k(s)}{s+k}\ .
\end{equation}
\begin{equation}
    \label{eq:main formula intro}
    \zeta (s) = \frac{1}{\Gamma(s)} \sum_0^\infty \frac{\alpha_k(s)}{s+k-1}\ .
\end{equation}
For positive integer $\lambda$:
\begin{eqnarray}
    \label{eq:3 intro}
    \zeta(s-\lambda)
    =
    \frac{1}{\Gamma(s)}
    \sum_{k=0}^\infty \alpha_k(s)
    \sum_{j=1}^\lambda
    (-1)^{\lambda+j} \frac{j! S(\lambda,j)}{s+k-j-1},
\end{eqnarray}
where $S(\lambda,j)$ are Stirling numbers of the second kind.
Finally,
\begin{equation}
    \label{eq:trigamma intro}
    \zeta(s+1)
    =
    \frac{1}{\Gamma(s)}
    \sum_{k=0}^\infty \alpha_k(s) \Psi_1(s+k),
\end{equation}
with $\Psi_1(s+k)$ the trigamma function.

At the time of writing~\cite{R}, the author also developed analogous formulas
for Dirichlet $L$-functions via the Hurwitz zeta function, defined in equation~\eqref{eq:hurwitz zeta},
as well as a variety
of additional formulas for the Riemann zeta function and the Gamma function,
and we report on these here. We collect our main formulas in the theorem below.
In equations~\eqref{eq:Gamma w intro}-~\eqref{eq:L with shift intro}, each stated
formula is valid for any $s \in \C$ that is not a pole of the lhs.

\begin{thm}

\label{thm:1}

Let $\Re{w}>0$. Then,
\begin{equation}
    \label{eq:Gamma w intro}
    \Gamma(s) = w^s \Gamma(w)
    \sum_{k=0}^\infty \frac{\alpha_k(s)\Gamma(s+k)}{\Gamma(s+k+w)}.
\end{equation}
Specializing to $w=N+1$, a positive integer:
\begin{equation}
    \label{eq:gamma N intro}
    \Gamma(s) = (N+1)^s N!
    \sum_{k=0}^\infty \frac{\alpha_k(s)}{(s+k)(s+k+1)\ldots(s+k+N)}.
\end{equation}
For $a>0$:
\begin{equation}
    \label{eq:hurwitz expansion intro}
    \zeta(s,a)
    = \frac{1}{\Gamma(s)}
    \sum_{k=0}^\infty \alpha_k(s) \frac{\Gamma(s+k-1)\Gamma(a)}{\Gamma(s+k+a-1)}.
\end{equation}
For non-negative integer $N$:
\begin{equation}
    \label{eq:zeta via hurwitz b intro}
    \zeta(s) =
    \sum_1^N n^{-s}
    + \frac{N!}{\Gamma(s)} \sum_{k=0}^\infty \frac{\alpha_k(s)}{ (s+k-1)\ldots(s+k-1+N)}.
\end{equation}
For $a>0$ and $\lambda$ a non-negative integer:
\begin{equation}
    \label{eq:hurwitz zeta with lambda shift intro}
    \zeta(s-\lambda,a)
    = \frac{1}{\Gamma(s)}
    \sum_{k=0}^\infty \alpha_k(s)
    \sum_{j=0}^\lambda
    \frac{\Gamma(s+k-j-1)\Gamma(a)}{\Gamma(s+k+a-j-1)} c_a(\lambda,j),
\end{equation}
where the coefficients $c_a(\lambda,j)$ are defined by~\eqref{eq:c_a}.

We also have, for positive integer $\Lambda$:
\begin{equation}
    \label{eq:desired formula intro}
    \sum_{\lambda=1}^\Lambda b_\lambda \zeta(s-\lambda)
    = \frac{1}{\Gamma(s)}
    \sum_{k=0}^\infty \frac{\alpha_k(s)}{(s+k-2)\ldots(s+k-\Lambda-1)},
\end{equation}
where the $b_\lambda$ are given by~\eqref{eq:b_lambda}.

Let $\chi$ be a non-trivial Dirichlet character for the modulus $q$. Then,
\begin{equation}
    \label{eq:L expansion intro}
    L(s,\chi)
    = \frac{1}{q^s\Gamma(s)}
    \sum_{k=0}^\infty
    \alpha_k(s) \Gamma(s+k-1)
    \sum_{m=1}^{q-1}
    \frac{\chi(m) \Gamma(m/q)}{\Gamma(s+k+m/q-1)}.
\end{equation}
Furthermore, for $\lambda$ a non-negative integer:
\begin{eqnarray}
    \label{eq:L with shift intro}
    L(s-\lambda,\chi)
    &=& \frac{1}{\Gamma(s) q^{s-\lambda}}
    \sum_{k=0}^\infty
    \alpha_k(s)
    \sum_{m=1}^{q-1}
    \chi(m) \Gamma(m/q) \notag \\
    &\times& \sum_{j=0}^\lambda
    \frac{\Gamma(s+k-j-1)}{\Gamma(s+k+m/q-j-1)} c_{m/q}(\lambda,j), \notag \\
\end{eqnarray}
where $c_{m/q}(\lambda,j)$ are defined by~\eqref{eq:c_a}

\end{thm}

As a consequence of these formulas, we have several identities for special
values. Two formulas for $L(s,\chi)$ at non-positive integers $s$
are summarized in the following theorem.

\begin{thm}
For non-negative integer $r$ and a non-trivial Dirichlet character $\chi$ for the modulus $q$:
\begin{equation}
    \label{eq:L -r intro}
    L(-r,\chi) = r! q^r
    \sum_{k=0}^{r} \frac{(-1)^{k-1}\alpha_k(-r)}{(r+1-k)!}
    \sum_{m=1}^{q-1} \chi(m) (m/q-1)\ldots(m/q+k-r-1).
\end{equation}
and, for non-negative integer $\lambda$:
\begin{eqnarray}
    \label{eq:L formula 2 intro}
    L(1-\lambda,\chi)
    &=& q^{\lambda-1}
    \sum_{j=0}^\lambda
    \left(
        \frac{(-1)^j}{j!}
        + \sum_{k=1}^j
        \frac{(-1)^{k-j}}{(j-k)!} \alpha_k'(1)
    \right) \notag \\
    &\times& \sum_{m=1}^{q-1}
    \chi(m) c_{m/q}(\lambda,j) (m/q-1)\ldots(m/q-j+k). \notag \\
\end{eqnarray}

\end{thm}

To compute $\alpha_k'(1)$ one can use either of the formulas from~\cite{R}:
\begin{equation}
    \label{eq:alpha prime}
    \alpha_{k+1}(1)' = \frac{1}{k+2}
    -\frac1{k+1} \sum_{j=1}^k
        \frac{j}{k-j+2} \alpha_j(1)'\ , \qquad k \geq 0.
\end{equation}
or
\begin{equation}
    \alpha_k'(1) = \frac{1}{k\ k!}
    \int_{0}^1 (x)_k dx,
\end{equation}
where
\begin{equation}
    (x)_k = x(x+1)\ldots(x+k-1)\ .
\end{equation}

We also derive several other formulas for special values. For example,
\begin{equation}
    \label{eq:gamma eqn intro}
    \gamma = \sum_{m=1}^N \frac{1}{m} -\log(N+1)
    -N! \sum_{k=1}^\infty \frac{\alpha_k(0)}{k(k+1)\ldots(k+N)},
\end{equation}
is shown, in Section~\ref{sec:euler gamma}, to follow from~\eqref{eq:gamma N intro}.
Equation~\eqref{eq:gamma eqn intro} is known and first due to Kluyver~\cite{Kl}.
In~\cite{R}, the author attributed his inspiration to Kenter's short note on
$\gamma$~\cite{Ke}, but the author has since discovered that Kluyver~\cite{Kl} essentially had
the same formulas as Kenter for $\gamma$, and also~\eqref{eq:gamma eqn intro}. The
above can also be viewed as a more precise form of the formula, often taken as
the definition of $\gamma$,
\begin{equation}
    \label{eq:gamma usual}
    \gamma = \lim_{N \to \infty} \sum_{m=1}^N \frac{1}{m} -\log(N+1).
\end{equation}

Finally, we derive several interesting formulas that generalize a formula of Hasse.

The alternating zeta function, or Dirichlet eta function, is defined by the Dirichlet series:
\begin{equation}
    \eta(s)
    := \sum_{n=1}^{\infty} \frac{(-1)^{n-1}}{n^{s}} = \zeta(s) (1-2^{1-s}).
\end{equation}

\begin{thm}
For all $s \in \C$ , and $\Re s_0 > -1$:
\begin{eqnarray}
    &&\eta(s+s_0) = \notag \\
    &&\frac{1}{\Gamma(s_0+1)}
    \sum_{m=0}^\infty (m+1)
    \left(
        \int_0^\infty \frac{x^{s_0}\exp(-x(m+1))}{(1+\exp(-x))^{m+2}} dx
    \right)
    \sum_{j=0}^m \frac{(-1)^j {m \choose j}}{(j+1)^s}.\notag \\
    \label{eq:hasse generalized}
\end{eqnarray}

Furthermore, for all $s,s_0 \in \C$:
\begin{eqnarray}
    \label{eq:zeta alt thm}
    &&\eta(s+s_0) = \notag \\
    &&
    \sum_{m=0}^\infty \frac{(-1)^m}{m!}
    \sum_{l=1}^{m+1} s(m+1,l) \eta(s_0+1-l)
    \sum_{j=0}^m \frac{(-1)^j {m \choose j}}{(j+1)^s}. \notag \\
\end{eqnarray}
where $s(m+1,l)$ are Stirling numbers of the first kind.
\end{thm}

Notice that~\eqref{eq:zeta alt thm} expresses $\eta(s)$, and hence $\zeta(s)$,
in terms of its values at any collection of points $s_0,s_0-1,s_0-2,\ldots$.

The special case $s_0=0$ in the above theorem simplifies.
Substituting $t=\exp(-x)$ into the integral of the first formula,
\begin{eqnarray}
    \label{eq:simplify s_0=0}
    &&(m+1)
    \int_0^\infty
    \frac{\exp(-x(m+1))}{(1+\exp(-x))^{m+2}} dx
    =
    (m+1) \int_0^1 \frac{t^m}{(1+t)^{m+2}} dt \notag \\
    &=& \frac{t^{m+1}}{(1+t)^{m+1}} \bigg|^1_0 =\frac{1}{2^{m+1}},
\end{eqnarray}
we get
\begin{equation}
    \label{eq:hasse}
    \eta(s) = \zeta(s)(1-2^{1-s}) =
    \sum_{m=0}^\infty \frac{1}{2^{m+1}}
    \sum_{j=0}^m \frac{(-1)^j {m \choose j}}{(j+1)^s},
\end{equation}
a formula first conjectured by Knopp and proved by Hasse~\cite{H}, later
rediscovered by Sondow~\cite{S}.

A similar formula holds for Dirichlet $L$-functions for a limited number of Dirichlet
characters. For any non-trivial $\chi$ mod $q\leq 5$ we have:
\begin{equation}
    L(s+s_0,\chi) =
    \sum_{m=0}^\infty \frac{(-1)^m}{m!}
    \sum_{l=1}^{m+1} s(m+1,l) L(s_0+1-l,\chi)
    \sum_{j=0}^m \frac{(-1)^j {m \choose j}}{(j+1)^s},
    \label{eq:L hasse}
\end{equation}
valid for all $s,s_0 \in \C$. The reason that we need to restrict here to
$q\leq 5$ has to do with the location of the $q$-th roots of unity relative to
the point 1, as described in Section~\ref{sec:dirichlet hasse}.

\subsection{A summation formula}
We give several other formulas related to the above. For example, we develop a
formula that holds for all non-trivial Dirichlet characters for any modulus. To
this end we use the following summation formula which we detail in
Section~\ref{sec:summation}.

Let $g,h$ be complex valued functions on the non-negative integers, say increasing
at most exponentially, i.e. assume there exists $c \geq 0$ such that $g(n),h(n)=O(\exp(cn))$.
Define the $m$-th finite differences recursively:
\begin{eqnarray}
    \Delta^0 h(j) &=&  h(j) \notag \\
    \Delta^m h(j) &=&  \Delta^{m-1} h(j+1) - \Delta^{m-1} h(j), \qquad m \geq 1.
\end{eqnarray}
Then
\begin{equation}
    \sum_{j=0}^m (-1)^j h(j+1) {m \choose j}
    = (-1)^m \Delta^m h(1).
    \label{eq:identity finite diff}
\end{equation}
Let
\begin{equation}
    G(z):=\sum_{n=1}^\infty g(n) z^{n-1}.
\end{equation}
We show that
\begin{equation}
    \sum_{n=1}^\infty g(n)h(n) z^{n-1}
    =
    \sum_{m=0}^\infty G^{(m)}(z) \Delta^{m} h(1) \frac{z^m}{m!},
    \label{eq:sum1}
\end{equation}
holds in some disc centred on $z=0$.

In some applications,
for example to the Hurwitz zeta function,
one might prefer to write the $n$-th term as $g(n)h(n) z^{n}$ and start the sum at $n=0$.
Thus, for functions $g,h$ on the non-negative integers, growing at most exponentially,
\begin{equation}
    \sum_{n=0}^\infty g(n)h(n) z^n
    =
    \sum_{m=0}^\infty \tilde{G}^{(m)}(z) \Delta^{m} h(0) \frac{z^m}{m!},
    \label{eq:sum2}
\end{equation}
holds on some disc centred on $z=0$,
where
\begin{equation}
    \tilde{G}(z):=\sum_{n=0}^\infty g(n) z^{n}.
\end{equation}

The analytic continuation of the rhs of~\eqref{eq:sum1} or~\eqref{eq:sum2} to a given
point $z$ depends on the location of the singularities of $G(w)$ or $\tilde{G}(w)$
in relation to the point $z$, and also on the rate of growth of $\Delta^m h(1)$ or
$\Delta^m h(0)$. In our applications, we substitute $z=1$ after considering
the analytic continuation of the function $G$ or $\tilde{G}$.

As a simple illustrative example, take $g(n)=1$ if $n=N$ and 0 otherwise, so that $\tilde{G}(z)=z^{N}$,
and $\tilde{G}^{(m)}(z)=0$ for all $m > N$. Substituting $z=1$ we get, as a special
case, Newton's forward difference formula:
\begin{equation}
    \label{eq:Newton}
    h(N) =  h(0) + N \Delta h(0) +\frac{N(N-1)}{2!} \Delta^2 h(0) +
    \frac{N(N-1)(N-2)}{3!} \Delta^3 h(0) +\ldots .
\end{equation}

As an application, we give the following identity.
\begin{thm}
Let $\lambda > 0$. If $\Re{s}>0$, then
\begin{equation}
    \label{eq:amore variant thm}
    \eta(s) = \exp(-\lambda) \sum_{m=0}^\infty \lambda^m
    \sum_{j=0}^m \frac{(-1)^j}{(j+1)^s \lambda^j (m-j)!},
\end{equation}
the sum over $m$ being uniformly convergent on compact subsets of $\Re{s}>0$.
Furthermore, as $m \to \infty$,
\begin{eqnarray}
    \sum_{j=0}^m \frac{(-1)^j}{(j+1)^s \lambda^j (m-j)!}
    \sim \frac{(-1)^m \exp(-\lambda)}{\lambda^m (m+1+\lambda)^s},
    \label{eq:asympt amore variant}
\end{eqnarray}
uniformly on compact subsets of $\Re{s}>0$.
\end{thm}

Our summation formula is then applied to the Hurwitz zeta function in Section~\ref{sec:summation hurwitz}
and to Dirichlet $L$-functions in Section~\ref{sec:summation dirichlet}.

For the Hurwitz zeta function, we actually consider
the alternating Hurwitz zeta defined by
\begin{equation}
    \zeta^*(s,a) :=
    \sum_{n=0}^{\infty} \frac{(-1)^{n}}{(n+a)^{s}},
\end{equation}
where $a>0$, and, in the above series, $\Re{s}>0$.
Our summation formula gives the following identities, as well as analytic continuation in
the stated regions below.
\begin{thm}
For $\Re{s_0}>-1$ and $s \in \C$,
\begin{eqnarray}
    &&\zeta^*(s+s_0,a) = \notag \\
    &&\frac{1}{\Gamma(s_0+1)}
    \sum_{m=0}^\infty 
    \left(
        \int_0^\infty \frac{x^{s_0}\exp(-x(m+a))(m+a+(a-1)\exp(-x))}{(1+\exp(-x))^{m+1}} dx
    \right) \notag \\
    && \qquad \qquad \qquad \times \sum_{j=0}^m \frac{(-1)^j {m \choose j}}{(j+a)^s}.
    \label{eq:another eq}
\end{eqnarray}
For all $s_0, s \in \C$,
\begin{equation}
    \zeta^*(s+s_0,a) =
    \sum_{m=0}^\infty \frac{(-1)^m}{m!}
    \sum_{l=1}^{m+1} s(m+1,l) \zeta^*(s_0+1-l,a)
    \sum_{j=0}^m \frac{(-1)^j {m \choose j}}{(j+a)^s}.
\end{equation}
Let $N+1$ be a positive integer.
Then, for all $s_0, s \in \C$,
\begin{eqnarray}
    &&\zeta^*(s+s_0,N+1) = \notag \\
    &&\sum_{m=0}^\infty \frac{(-1)^m}{m!}
    \sum_{l=1}^{m+1} s(m+1,l) \zeta^*(s_0+1-l,N+1)
    \sum_{j=0}^m \frac{(-1)^j {m \choose j}}{(j+1)^s}. \notag \\
\end{eqnarray}
\end{thm}
Notice that the second formula above is {\it not} a specialization of the first.
The difference between the above two formulas
is in the factor $j+1$ vs $j+a$ in the sum over $j$.

For Dirichlet $L$-functions we have the following theorem, valid for any
non-trivial Dirichlet character for any modulus $q$.
\begin{thm}
For $\Re{s_0}>-1$, and $s \in \C$,
\begin{equation}
    L(s+s_0,\chi) =
    \sum_{m=0}^\infty (m+1)
    \left(
        \int_0^\infty \frac{x^{s_0}\exp(-x(m+1))}{(1+\exp(-x))^{m+2}} dx
    \right)
    \sum_{j=0}^m \frac{\chi(j+1) {m \choose j}}{(j+1)^s},
\end{equation}
and, for all $s_0,s \in \C$,
\begin{equation}
    L(s+s_0,\chi) =
    \sum_{m=0}^\infty \frac{(-1)^m}{m!}
    \sum_{l=1}^{m+1} s(m+1,l) \eta(s_0+1-l)
    \sum_{j=0}^m \frac{\chi(j+1) {m \choose j}}{(j+1)^s}.
\end{equation}
\end{thm}
Notice that these expansions express $L(s,\chi)$ in terms of the eta function
$\eta(s)=\zeta(s)(1-2^{1-s})$!

As a special case, set $s_0=0$ to get
\begin{equation}
    L(s,\chi) =
    \sum_{m=0}^\infty \frac{1}{2^{m+1}}
    \sum_{j=0}^m \frac{\chi(j+1) {m \choose j}}{(j+1)^s},
    \label{eq:L method 2 c}
\end{equation}
valid for any non-trivial Dirichlet character $\chi$, and
all $s \in \C$.

Formulas~\eqref{eq:hasse generalized}-~\eqref{eq:L hasse} and
~\eqref{eq:another eq}-~\eqref{eq:L method 2 c} converge
uniformly for $s_0,s$ in compact subsets of the regions in which the formulas
,are claimed to be true. As part of their derivation we give bounds for the sums over $j$
that appear in these formulas.

Let $\tau(\chi)$ be the Gauss sum
\begin{equation}
    \label{eq:gauss sum}
    \tau(\chi) = \sum_{j=0}^{q-1} \chi(j) e(j/q),
\end{equation}
where
\begin{equation}
    e(t):=\exp(2\pi i t).
\end{equation}

We have the following estimates.

\begin{thm}
Let $a>0$. Then, as $m \to \infty$,
\begin{equation}
    \label{eq:j sum asympt}
    \sum_{j=0}^m \frac{(-1)^j}{(j+a)^s} {m \choose j} \sim
    \frac{\log(m)^{s-1}\Gamma(a)}{m^a\Gamma(s)},
\end{equation}
uniformly for $s$ in compact subsets of $\C$. In the event that
$s \in \Z$ and $s\leq 0$, then we interpret the $\sim$ to mean equality when
$m>|s|$.

Let $\chi$ be a non-trivial Dirichlet character for the modulus $q$. Then as $m/q^2 \to \infty$
\begin{eqnarray}
    &&\sum_{j=0}^m \chi(j+1) {m \choose j} \sim \notag \\
    &&\frac{\tau(\chi)}{q} \left(e(-1/q)(1+e(-1/q))^m+ \chi(-1)e(1/q) (1+e(1/q))^m\right)\notag \\
    &&\leq \frac{2}{q^{1/2}} C_q^m,
    \label{eq:chi sum bound thm}
\end{eqnarray}
and
\begin{equation}
    C_q=|1+e(1/q)|<2.
    \label{eq:C_q b}
\end{equation}
Furthermore,
\begin{eqnarray}
    &&\sum_{j=0}^m \chi(j+1) {m \choose j} \notag \\
    &&\sim \frac{\tau(\chi)}{q} 2^{m+1} e^{-\frac{m\pi^2}{2q^2}}
    \times
    \begin{cases}
         \cos({\frac{\pi}{q}(m+2)}), &\text{if $\chi(-1)=1$,} \\
        -i\sin({\frac{\pi}{q}(m+2)}),&\text{if $\chi(-1)=-1$,} \notag \\
    \end{cases} \\
    \label{eq:chi sum asympt}
\end{eqnarray}
as $m/q^2 \to \infty$ and $m/q^4 \to 0$. In~\eqref{eq:chi sum bound thm}
and~\eqref{eq:chi sum asympt}, we interpret $\sim$ to mean equality whenever
the rhs vanishes.

We also have, for $q$ fixed and as $m \to \infty$,
\begin{eqnarray}
    &&\sum_{j=0}^m \frac{\chi(j+1)}{(j+1)^s} {m \choose j} \ll m^M C_q^m,
    \label{eq:chi sum bound any s}
\end{eqnarray}
uniformly for $s$ on compact subsets of $\Re{s}>-M$, where $M$ is a non-negative integer,
and with the implied constant in the $\ll$ also depending on $q$ and $M$.

\end{thm}
In light of identity~\eqref{eq:identity finite diff},
these estimates say that the $m$-th finite differences of $1/(j+a)^s$ and of $(-1)^j
\chi(j+1)/(j+1)^s$ are exponentially smaller than the trivial bound obtained by
using $\sum_0^m {m \choose j}=2^m$.
In the case of Dirichlet characters,
the extra power of $-1$ is crucially needed to get
such a small bound.

We remark that Coffey~\cite{C}, building on the author's paper~\cite{R}, has
separately arrived at formulas~\eqref{eq:Gamma w intro} and ~\eqref{eq:hurwitz expansion intro}
of Theorem~\ref{thm:1}. His application, however, is to the polygamma function
and Stieltjes constants, whereas our focus is on formulas for $L$-functions and
analytic aspects. Furthermore, his derivation is formal in that it does not
contain any discussion on the convergence, or rate of convergence, of the sums
involved in these formulas, and merely asserts analytic continuation.

Coffey uses, in his paper, the expansion
\begin{equation}
    \label{eq:Norlund 1}
    \left(\frac{\ln(x+1)}{x}\right)^z=z \sum_{k=0}^\infty \frac{x^k}{k!}
    \frac{B_k^{(z+k)}}{(z+k)}, \qquad |x| <1,
\end{equation}
where the coefficients $B_k^{(z+k)}$ are generalized Bernoulli numbers given by
the generating function
\begin{equation}
    \left({z \over {e^z-1}}\right)^\alpha =\sum_{n=0}^\infty
    B_n^{(\alpha)} \frac{z^n}{n!}, \qquad |z| < 2\pi,
    \label{eq:Norlund 2}
\end{equation}
and remarks that
\begin{equation}
    \alpha_k(s)=(s-1)\frac{(-1)^k}{k!}\frac{B_k^{(s+k-1)}}{(s+k-1)}.
    \label{eq:alpha vs B}
\end{equation}
We regard the generating function~\eqref{eq:1} and the coefficients $\alpha_k(s)$,
to be more fundamental than the generating functions~\eqref{eq:Norlund 2} and~\eqref{eq:Norlund 1}
to the problem at hand.

In this paper we focus our attention on deriving formulas for $L$-functions of
degree 1. It seems, however, that some of our techniques can be applied to
higher degree $L$-functions. We plan to explore higher degree $L$-functions,
and to consider computational aspects of these formulas and variants, in a
subsequent paper.

\section{Expansion for the Hurwitz zeta function}

The Hurwitz zeta function is defined, for $\Re{s}>1$ and $a>0$, by
\begin{equation}
    \label{eq:hurwitz zeta}
    \zeta(s,a) = \sum_{n=0}^\infty (n+a)^{-s}.
\end{equation}
We have
\begin{eqnarray}
    \label{eq:hurwitz zeta manipulate} \notag
    \Gamma(s) \zeta(s,a) &=& \int_0^\infty x^{s-1} \sum_{n=0}^\infty \exp(-(n+a)x) dx \\
    &=& \int_0^\infty \frac{x^{s-1} \exp(-ax)}{1-\exp(-x)} dx.
\end{eqnarray}
Substituting $t=1-\exp(-x)$, the above becomes, on applying~\eqref{eq:1},
\begin{eqnarray}
    \label{eq:hurtwitz manipulate} \notag
    \Gamma(s) \zeta(s,a) &=& \int_0^1 (-\log(1-t))^{s-1} (1-t)^{a-1} \frac{dt}{t} \\
    &=& \sum_{k=0}^\infty \alpha_k(s) \int_0^1 t^{s+k-2}(1-t)^{a-1} dt.
\end{eqnarray}
The rearranging of integration and summation over $k$ is
justified by bound~\eqref{eq:alpha bound} which shows that, for given $s$,
the sum over $k$ converges uniformly for $t$ in any closed subset of $(0,1)$.
Now
\begin{equation}
    \label{eq:beta}
    \int_0^1 t^{s+k-2}(1-t)^{a-1} = \beta(s+k-1,a) = \frac{\Gamma(s+k-1)\Gamma(a)}{\Gamma(s+k+a-1)}.
\end{equation}
Therefore,
\begin{equation}
    \label{eq:hurwitz expansion}
    \zeta(s,a)
    = \frac{1}{\Gamma(s)}
    \sum_{k=0}^\infty \alpha_k(s) \frac{\Gamma(s+k-1)\Gamma(a)}{\Gamma(s+k+a-1)}.
\end{equation}
Note that~\eqref{eq:hurwitz expansion} combined with the bound~\eqref{eq:alpha bound}
for $|\alpha_k(s)|$ also provides the meromorphic continuation in $s$
of $\zeta(s,a)$ to all of $\C$, with the only pole being simple at $s=1$,
because the sum converges uniformly for $s$
in compact sets away from its poles. For the latter fact, we also need
\begin{equation}
    \label{eq:k ratio}
    \frac{\Gamma(s+k-1)}{\Gamma(s+k+a-1)} \sim \frac{1}{k^a},
\end{equation}
as $k \to \infty$ which follows from Gauss' formula~\eqref{eq:Gauss limit},
also expressible as
\begin{equation}
    \label{eq:Gauss limit b}
    \Gamma(s+N+1) \sim N^s N!.
\end{equation}
Apply the above with $N$ replaced by $k$, and $s$ replaced, in the numerator
of~\eqref{eq:k ratio}, by $s-2$, and, in the denominator, by $s+a-2$.

The poles of the numerator of the $k$-th summand in~\eqref{eq:hurwitz expansion} occur
at $s=-k+1,-k,-k-1,-k-2,\ldots$. All of these, except for the pole at $s=1$ when $k=0$,
are cancelled by the zeros of $1/\Gamma(s)$ at $s=0,-1,-2,\ldots$, thus the only pole
occurs at $s=1$.

If we let, above, $a=1$ and use $\Gamma(k+s-1)/\Gamma(s+k) = 1/(s+k-1)$ we
recover~\eqref{eq:main formula intro}. And, if we take $a=N+1$ to be
a positive integer, we get
\begin{equation}
    \label{eq:zeta via hurwitz}
    \zeta(s) - \sum_1^N  n^{-s}
    = \frac{N!}{\Gamma(s)} \sum_{k=0}^\infty \frac{\alpha_k(s)}{ (s+k-1)\ldots(s+k-1+N)},
\end{equation}
i.e.
\begin{equation}
    \label{eq:zeta via hurwitz b}
    \zeta(s) =
    \sum_1^N n^{-s}
    + \frac{N!}{\Gamma(s)} \sum_{k=0}^\infty \frac{\alpha_k(s)}{ (s+k-1)\ldots(s+k-1+N)}.
\end{equation}
We can increase the rate of convergence of this sum by selecting $N$ to be
larger. For example, in Table~\ref{tab:1} we compare the precision achieved
with $s=3$, $N=1$, $N=5$, and $N=20$, $N=100$, and $k\leq K$, for various values of $K$.
More specifically, we truncate the sum at $k\leq K$ and denote the remainder by :
\begin{equation}
    \label{eq:R}
    R(K,N,s) := \frac{N!}{\Gamma(s)} \sum_{k=K+1}^\infty \frac{\alpha_k(s)}{ (s+k-1)\ldots(s+k-1+N)}.
\end{equation}
In our table, we depict the relative remainder terms
\begin{equation}
    \label{eq:R tilde}
    \tilde{R}(K,N,s):= R(K,N,s)/\zeta(s),
\end{equation}
for $s=3$.

\begin{table}[h!]
\centerline{
\begin{tabular}{|c|c|c|c|c|} \hline
    $K$ & $\tilde{R}(K,1)$ & $\tilde{R}(K,5)$ & $\tilde{R}(K,20)$ & $\tilde{R}(K,100)$ \\ \hline
0 & 0.09876701304 & 0.003735350802 & 8.885612538e-05 & 8.054816065e-07 \\ \hline
1 & 0.06410420585 & 0.001259436003 & 1.056632936e-05 & 2.148135338e-08 \\ \hline
2 & 0.04503966189 & 0.0005029064811 & 1.595623563e-06 & 7.505774529e-10 \\ \hline
4 & 0.02594760777 & 0.0001135959722 & 6.161196116e-08 & 1.697998088e-12 \\ \hline
8 & 0.01207251038 & 1.310916493e-05 & 3.577149365e-10 & 4.694895745e-17 \\ \hline
16 & 0.004670251882 & 7.985644458e-07 & 2.203486644e-13 & 1.938746948e-24 \\ \hline
32 & 0.001580885488 & 2.957659659e-08 & 1.468535646e-17 & 1.343384486e-35 \\ \hline
64 & 0.0004908542321 & 7.867242541e-10 & 1.559347771e-22 & 3.318995397e-51 \\ \hline
128 & 0.0001445094596 & 1.715803026e-11 & 4.578210941e-28 & 2.669903582e-71 \\ \hline
256 & 4.117178337e-05 & 3.342699708e-13 & 6.093017067e-34 & 3.364327536e-95 \\ \hline
512 & 1.148739118e-05 & 6.112998825e-15 & 5.183860931e-40 & 7.923516397e-122 \\ \hline
1024 & 3.160289401e-06 & 1.078325101e-16 & 3.464372739e-46 & 3.670183638e-150 \\ \hline
2048 & 8.606974887e-07 & 1.861800819e-18 & 2.037483910e-52 & 1.820228675e-179 \\ \hline
\end{tabular}
}
\caption
{A comparison of the remainder terms $\tilde{R}(K,N,3)$ for $N=1,5,20,100$ and
various values of $K$.}
\label{tab:1}
\end{table}

Out of curiousity, we compare with the traditional way of computing $\zeta(s)$
using Euler-MacLaurin summation:
\begin{eqnarray}
    \label{eq:zeta euler mac}
    \zeta(s) &=&
    \sum_1^N n^{-s} +
    \frac{N^{1-s}}{s-1}
    +\sum_1^K {s+k-2 \choose k-1} \frac{B_k}{k} N^{-s-k+1} \notag \\
    &-&{s+K-1 \choose K} \int_N^\infty B_K(\{t\}) t^{-s-K} dt,
\end{eqnarray}
where $B_k, B_k(x)$ denotes the Bernoulli numbers and polynomials \cite{R2}.
and $B_k:=B_k(0)$. Table~\ref{tab:2} depicts, for $s=3$ and the same values of
$N$ and $K$ as in Table~\ref{tab:1}, the relative remainder term
$\tilde{R}_2(K,N,s):=R_2(K,N,s)/\zeta(s)$, where
\begin{equation}
    \label{eq:R2}
    R_2(K,N,s):=
    -{s+K-1 \choose K} \int_N^\infty B_K(\{t\}) t^{-s-K} dt.
\end{equation}

\begin{table}[h!]
\centerline{
\begin{tabular}{|c|c|c|c|c|} \hline
    $K$ & $\tilde{R}_2(K,1)$ & $\tilde{R}_2(K,5)$ & $\tilde{R}_2(K,20)$ & $\tilde{R}_2(K,100)$ \\ \hline
0 & -0.2478610589 & -0.002999137452 & -5.069543603e-05 & 4.138739872e-07 \\ \hline
1 & 0.1680926274 & 0.0003284920385 & 1.298774753e-06 & -2.079699113e-09 \\ \hline
2 & -0.03988421573 & -4.270910495e-06 & -1.080516795e-09 & 6.931868307e-14 \\ \hline
4 & 0.02944139866 & 1.659288259e-07 & 2.695929470e-12 & -6.931313924e-18 \\ \hline
8 & 0.08490189016 & 1.233350305e-09 & 8.380153455e-17 & -3.464912761e-25 \\ \hline
16 & 43.97791945 & 3.248852315e-12 & 4.036386548e-24 & -4.339822180e-38 \\ \hline
32 & 2.000506814e+11 & 1.894001979e-13 & 8.399102960e-35 & -6.233469067e-60 \\ \hline
64 & 5.607820702e+39 & 3.383888667e-07 & 1.661297238e-47 & -6.266500089e-95 \\ \hline
128 & 2.810377072e+115 & 3.618530655e+24 & 1.065733432e-54 & -1.163548513e-146 \\ \hline
256 & 8.493688329e+304 & 3.878203997e+124 & 1.480620677e-31 & -1.227926982e-212 \\ \hline
512 & 3.186217020e+760 & 1.702666495e+401 & 5.605030930e+91 & -1.276702408e-268 \\ \hline
1024 & 2.130486030e+1825 & 1.530602415e+1108 & 2.919002699e+490 & -1.548683921e-227 \\ \hline
2048 & 6.062592955e+4262 & 7.835225059e+2829 & 4.672401590e+1595 & -5.541412864e+162 \\ \hline
\end{tabular}
}
\caption
{A comparison of the remainder terms $\tilde{R}_2(K,N,3)$ for $N=1,5,20,100$ and
various values of $K$.}
\label{tab:2}
\end{table}

Formulas~\eqref{eq:zeta via hurwitz b} and~\eqref{eq:zeta euler mac} are
different in a number of ways. First, our sum over $k$ converges and truncating the
series after $K$ terms gives a better approximation as $K \to \infty$, while
the sum in the Euler-MacLaurin formula provides a divergent
asymptotic expansion for $\zeta(s)$.
Second, our formulas give the meromorphic continuation of $\zeta(s)$ and of
$\zeta(s,a)$ to all of $\C$ with a simple pole with residue 1 at $s=1$. On the other hand,
the Euler-MacLaurin formula only gives, for given $K$, the meromorphic continuation
of $\zeta(s)$ up to $\Re{s}>-K+1$.

One can show (see 2.2.5 of~\cite{R2}), that the remainder term in the
Euler-MacLaurin formula satisfies, for $\sigma:=\Re{s} > -K+1$,
\begin{equation}
    \label{eq: R2 estimate}
    |R_2(K,N,s)| \leq
    \frac{2\zeta(K)}{N^{\sigma-1}}
    \frac{|s+K-1|}{\sigma+K-1}
    \frac{|\Gamma(s+K-1)|}{|\Gamma(s)| (2\pi N)^K}.
\end{equation}
On the other hand, writing the $k$-th term in~\eqref{eq:R} as
\begin{equation}
    \label{eq:zeta term}
    \frac{\alpha_k(s) \Gamma(s+k-1) \Gamma(N+1)}{\Gamma(s) \Gamma(s+k+N)}
\end{equation}
we have, for given $s$ and $K$, that
\begin{equation}
    \label{eq:R N limit}
    |R(K,N,s)| = O(N^{-\sigma-K}),
\end{equation}
with the implied constant in the big-O depending on $s$ and $K$. Here we
have used the asymptotic formula
\begin{equation}
    \label{eq:ratio limit}
    \frac{\Gamma(N+1)}{\Gamma(s+k+N)} \sim \frac{1}{N^{s+k-1}}
\end{equation}
as $N\to \infty$, which can be proven using Gauss' formula~\eqref{eq:Gauss limit}
with $s$ replaced by $s+k-1$ or using Stirling's formula. Therefore, in
the $N$-aspect, formula~\eqref{eq:zeta via hurwitz} is roughly comparable to
the Euler-MacLaurin formula, though with different dependence on $K$ and $s$.
The main difference between the two is in the extra $(2\pi)^{-K}$ in
\eqref{eq: R2 estimate}, and the $\alpha_k(s)$, which grows exponentially in $|s|$,
in each term of $R(K,N,s)$.

Next we examine the behaviour of $R(K,N,s)$ in the $K$-aspect.
The $k$-th term of~\eqref{eq:R} is
\begin{equation}
    \label{eq:kth term}
    \frac{\alpha_k(s) \Gamma(s+k-1) \Gamma(N+1)}{\Gamma(s) \Gamma(s+k+N)}.
\end{equation}
By~\eqref{eq:alpha bound},
this equals, for fixed $s$ and $N$,
\begin{equation}
    O(k^{-N-2+\epsilon}),
\end{equation}
for any $\epsilon >0$, because $\alpha_k(s) = O(k^{-1+\epsilon})$ and
$\Gamma(s+k-1)/\Gamma(s+k+N) = O(k^{-N-1})$. Summing over $k\geq K+1$ thus gives
\begin{equation}
    \label{eq:K estimate}
    R(K,N,s) = O(K^{-N-1+\epsilon}),
\end{equation}
with the implied constant depending on $s$ and $N$.

For comparison, we examine~\eqref{eq:zeta euler mac}, the Euler-MacLaurin
summation formula, in the $K$-aspect. In practice,~\eqref{eq:zeta euler mac} is
used to compute $\zeta(s)$
by summing the first row of the formula and dropping the term in the second row at
a cost that is estimated.
Now, for even $k$, $|B_k|/k! = 2\zeta(k)/(2\pi)^k$. Thus, for even $k$,
a typical term equals, in absolute value,
\begin{equation}
    \frac{\zeta(k)}{\pi N^{\sigma}} \prod_{j=0}^{k-2} \frac{|s+j|}{2\pi N}.
\end{equation}
For sufficiently large $j$, the numerator in the product can be as large as we wish compared to the
denominator,
i.e. for given $s$ and $N$, the terms in the sum over $k$ in~\eqref{eq:zeta euler mac} get arbitrarily
large as $k$ grows.

\section{A related expansion for the Gamma function}

Let $\Re{w}, \Re{s}>0$ and consider
\begin{eqnarray}
    \label{eq:gamma}
    \Gamma(s) w^{-s} &=& \int_0^\infty x^{s-1}e^{-wx} dx \notag \\
    &=& \int_0^1 (-\log (1-t)/t)^{s-1} t^{s-1}(1-t)^{w-1} dt \notag \\
    &=& \sum_{k=0}^\infty \alpha_k(s) \int_0^1 t^{s+k-1}(1-t)^{w-1} dt \notag \\
    &=& \sum_{k=0}^\infty \alpha_k(s) \frac{\Gamma(s+k) \Gamma(w)}{\Gamma(s+k+w)}.
\end{eqnarray}
Therefore,
\begin{equation}
    \label{eq:Gamma w}
    \Gamma(s) = w^s \Gamma(w)
    \sum_{k=0}^\infty \frac{\alpha_k(s)\Gamma(s+k)}{\Gamma(s+k+w)}.
\end{equation}
By~\eqref{eq:alpha bound}, this series converges uniformly for all $s$ in a
given compact set away from its poles and therefore provides the meromorphic
continuation of $\Gamma(s)$. Specializing to $w=N+1$, where $N$ is a
non-negative integer, we get
\begin{equation}
    \label{eq:gamma N}
    \Gamma(s) = (N+1)^s N!
    \sum_{k=0}^\infty \frac{\alpha_k(s)}{(s+k)(s+k+1)\ldots(s+k+N)}.
\end{equation}
Notice that, as $N\to \infty$, the $k=0$ term dominates. We can see this by
writing the $k$-th term of the sum as:
\begin{equation}
    \label{eq:gauss k term}
    \frac{\alpha_k(s)}{(s+k)}
    \frac{1}{s(s+1)\ldots(s+N)}
    \frac{s(s+1)\ldots(s+k-1)}{(s+N+1)\ldots(s+N+k)}.
\end{equation}
Now, $s/(s+N+1) = O(1/N)$, with the implied constant depending on $s$.
Furthermore, if $N$ is sufficiently large in comparison to $|\Re{s}|$,
then $|(s+j)/(s+N+j+1)| < 1$ for all positive integers $j$. Thus the terms
with $k \geq 1$ contribute
\begin{equation}
    \ll
    \frac{(N+1)^s N!}{s(s+1)\ldots(s+N)}
    \frac{1}{N} \sum_{k=1}^\infty
    \frac{|\alpha_k(s)|}{|s+k|}
\end{equation}
to~\eqref{eq:gamma N}.
Bound~\eqref{eq:alpha bound} shows that the sum over $k$ converges, hence
the above is $O(1/N)$ times the $k=0$ term. We thus
get Gauss' formula for $\Gamma(s)$:
\begin{equation}
    \label{eq:Gauss limit}
    \Gamma(s) = \lim_{N\to\infty} \frac{N^s N!}{s(s+1)\ldots(s+N)}.
\end{equation}
Here we have also used $(N+1)^s \sim N^s$ as $N\to \infty$.
Thus, the $k=0$ term in formula~\eqref{eq:gamma N} connects Euler's
formula for the Gamma function, i.e. our starting point, to that of Gauss.

\section{Formula for Euler's constant}
\label{sec:euler gamma}

In equation~\eqref{eq:gamma N}, subtract $1/s$ from both sides and then let $s\to 0$. On the lhs, we
get $-\gamma$, i.e. the negative of Euler's constant
(here, we can define $-\gamma$ to be the constant term in the Laurent series about $s=0$ of $\Gamma(s)$).
On the rhs, the $1/s$  will cancel with a portion of the $k=0$ term. More precisely, expand $(N+1)^s$
as:
\begin{equation}
    \label{eq:term1}
    (N+1)^s = 1 +\log(N+1) s + O(s^2).
\end{equation}
Furthermore, the $k=0$ term, with the $N!$ but not the $(N+1)^s$, equals
\begin{equation}
    \label{eq:k=0 term}
    \frac{N!}{s(s+1)\ldots(s+N)} = \sum_{m=0}^N \frac{(-1)^m}{s+m} {N \choose m},
\end{equation}
the latter from the partial fraction expansion~\eqref{eq:parfrac}. Therefore, multiplying~\eqref{eq:term1},
with~\eqref{eq:gamma N}, using the above for $k=0$, collecting terms, and taking the limit as $s \to 0$,
we have
\begin{equation}
    \gamma = \sum_{m=1}^N \frac{(-1)^{m+1}}{m} {N \choose m} -\log(N+1)
    -N! \sum_{k=1}^\infty \frac{\alpha_k(0)}{k(k+1)\ldots(k+N)}.
\end{equation}
Now, by equation~\eqref{eq:k=0 term}, the first sum on the rhs above equals
\begin{eqnarray}
    \left(
        \frac{1}{s}-
        \frac{N!}{s(s+1)\ldots(s+N)}
    \right)_{s=0}
    &=& \frac{(s+1)\ldots(s+N)-N!}{s(s+1)\ldots(s+N)} \bigg|_{s=0} \notag \\
    &=& \sum_{m=1}^N \frac{1}{m},
\end{eqnarray}
the last step by examining the coefficient of $s$ of the polynomial in the
numerator.
Hence
\begin{equation}
    \label{eq:gamma eqn}
    \gamma = \sum_{m=1}^N \frac{1}{m} -\log(N+1)
    -N! \sum_{k=1}^\infty \frac{\alpha_k(0)}{k(k+1)\ldots(k+N)}.
\end{equation}

\section{Hurwitz zeta function shifted}

Next we describe a formula for $\zeta(s-\lambda,a)$ where $\lambda$ is
a positive integer and $a>0$. Begin with
\begin{eqnarray}
    \label{eq:hurwitz zeta with shift}
    \Gamma(s) \zeta(s-\lambda,a)
    &=& \int_0^\infty x^{s-1} \sum_{n=0}^\infty (n+a)^\lambda \exp(-(n+a)x) dx \notag \\
    &=& \int_0^1 (-\log(1-t))^{s-1} \sum_{n=0}^\infty (n+a)^\lambda (1-t)^{n+a-1} dt. \notag \\
\end{eqnarray}
We can express
\begin{equation}
    \sum_{n=0}^\infty (n+a)^\lambda (1-t)^{n+a-1}
\end{equation}
as a rational function in $t$ by starting with the $\lambda=0$ case,
$\sum_{n=0}^\infty (1-t)^{n+a-1} = (1-t)^{a-1}/t$, and repeatedly
multiplying by $1-t$ and applying $-d/dt$. We can thus prove, inductively, that
\begin{equation}
    \label{eq:sum in closed form}
    \sum_{n=0}^\infty (n+a)^\lambda (1-t)^{n+a-1}
    = (1-t)^{a-1} \sum_{j=0}^\lambda \frac{c_a(\lambda,j)}{t^{j+1}},
\end{equation}
where $c_a(0,0)=1$, $c_a(0,1)=0$, and
\begin{equation}
    \label{eq:c_a}
    c_a(\lambda+1,j) = (a-j-1) c_a(\lambda,j) + j c_a(\lambda,j-1),
\end{equation}
for $\lambda \geq 1$, $0 \leq j \leq \lambda+1$ (for the recursion, we also set
$c_a(\lambda,\lambda+1)=0$. We need not worry about defining $c_a(\lambda,j-1)$
when $j=0$ because of the factor of $j$ that appears in front of the
$c_a(\lambda,j-1)$). As usual, we write
$(-\log(1-t))^{s-1} = \sum_{k=0}^\infty \alpha_k(s) t^{k+s-1}$,
and get
\begin{equation}
    \label{eq:hurwitz zeta with lambda shift}
    \zeta(s-\lambda,a)
    =\frac{1}{\Gamma(s)}
    \sum_{k=0}^\infty \alpha_k(s)
    \sum_{j=0}^\lambda
    \frac{\Gamma(s+k-j-1)\Gamma(a)}{\Gamma(s+k+a-j-1)} c_a(\lambda,j).
\end{equation}
We list a few $c_a(\lambda,j)$ in Table~\ref{tab:3}.
When $a=1$ we have
\begin{equation}
    c_1(\lambda,j) = (-1)^{\lambda+j} j! S(\lambda,j)
\end{equation}
and~\eqref{eq:hurwitz zeta with lambda shift} can be seen as a hybrid
of~\eqref{eq:3 intro} and~\eqref{eq:hurwitz expansion}.

\begin{table}[h!]
\centerline{
\begin{tabular}{|c|c|c|} \hline
    $\lambda$ & $j$ & $c_a(\lambda,j)$ \\ \hline
    0 & 0 & $1$ \\ \hline
    1 & 0 & $a-1$ \\
    1 & 1 & $1$ \\ \hline
    2 & 0 & $(a-1)^2$ \\
    2 & 1 & $2a-3$ \\
    2 & 2 & $2$ \\ \hline
    3 & 0 & $(a-1)^3$ \\
    3 & 1 & $3a^2-9a+7$ \\
    3 & 2 & $6a-12$ \\
    3 & 3 & $6$ \\ \hline
    4 & 0 & $(a-1)^4$ \\
    4 & 1 & $4a^3-18a^2+28a-15$ \\
    4 & 2 & $12a^2-48a+50$ \\
    4 & 3 & $24a-60$ \\
    4 & 4 & $24$ \\ \hline
    5 & 0 & $(a-1)^5$ \\
    5 & 1 & $5a^4-30a^3+70a^2-75a+31$ \\
    5 & 2 & $20a^3-120a^2+250a-180$ \\
    5 & 3 & $60a^2-300a+390$ \\
    5 & 4 & $120a-360$ \\
    5 & 5 & $120$ \\ \hline
    6 & 0 & $(a-1)^6$ \\
    6 & 1 & $6a^5-45a^4+140a^3-225a^2+186a-63$ \\
    6 & 2 & $30a^4-240a^3+750a^2-1080a+602$ \\
    6 & 3 & $120a^3-900a^2+2340a-2100$ \\
    6 & 4 & $360a^2-2160a+3360$ \\
    6 & 5 & $720a-2520$ \\
    6 & 6 & $720$ \\ \hline
    7 & 0 & $(a-1)^7$ \\
    7 & 1 & $7a^6-63a^5+245a^4-525a^3+651a^2-441a+127$ \\
    7 & 2 & $42a^5-420a^4+1750a^3-3780a^2+4214a-1932$ \\
    7 & 3 & $210a^4-2100a^3+8190a^2-14700a+10206$ \\
    7 & 4 & $840a^3-7560a^2+23520a-25200$ \\
    7 & 5 & $2520a^2-17640a+31920$ \\
    7 & 6 & $5040a-20160$ \\
    7 & 7 & $5040$ \\ \hline
\end{tabular}
}
\caption
{A table of $c_a(\lambda,j)$.}
\label{tab:3}
\end{table}

\section{A linear combination of $\zeta$ that converges quickly}

Let $\Lambda$ be a positive integer. In this section we derive a formula that
expresses a linear combination of $\zeta(s-1),\ldots,\zeta(s-\Lambda)$ in a
series of a nature similar to~\eqref{eq:zeta via hurwitz b}.

We will exploit the fact that the same coefficients $\alpha_k(s)$ appear in~\eqref{eq:3 intro}
independent of $\lambda$, and that the denominators in~\eqref{eq:3 intro} are
particularly simple. By taking linear combinations of $\Gamma (s)\zeta
(s-\lambda)$ we can develop a sum whose terms converge more rapidly. We do so
by using the partial fraction expansion:
\begin{equation}
    \label{eq:parfrac}
    \frac{1}{z(z-1)\ldots(z-m)} = \frac{a_0}{z} + \frac{a_1}{z-1} + \ldots \frac{a_m}{z-m},
\end{equation}
where
\begin{equation}
    \label{eq:a_l}
    a_l = \frac{(-1)^{m-l}}{l!(m-l)!} = \frac{(-1)^{m-l}}{m!} {m \choose l}.
\end{equation}
The formula for $a_l$ can be derived by considering the residue of both sides of~\eqref{eq:parfrac}
at $z=l$.

Therefore, letting $\Lambda$ be a positive integer, we wish to find a linear combination:
\begin{equation}
    \label{eq:linear comb}
    \sum_{\lambda=1}^\Lambda
        b_\lambda \zeta(s-\lambda)
    =
    \frac{1}{\Gamma(s)}
    \sum_{k=0}^\infty \alpha_k(s)
    \sum_{\lambda=1}^\Lambda
    \sum_{j=1}^\lambda
    (-1)^{\lambda+j} \frac{j! S(\lambda,j)}{s+k-j-1}
    b_\lambda
\end{equation}
such that
\begin{equation}
    \label{eq:goal}
    \sum_{\lambda=j}^\Lambda
    (-1)^{\lambda+j} S(\lambda,j) j! b_\lambda = \frac{(-1)^{\Lambda-j}}{(\Lambda-1)!} {\Lambda-1 \choose j-1},
\end{equation}
for each $1 \leq j \leq \Lambda$.
Doing so would then result, by~\eqref{eq:parfrac}, in the formula:
\begin{equation}
    \label{eq:desired formula}
    \sum_{\lambda=1}^\Lambda
        b_\lambda \zeta(s-\lambda)
    =
    \frac{1}{\Gamma(s)}
    \sum_{k=0}^\infty \frac{\alpha_k(s)}{(s+k-2)\ldots(s+k-\Lambda-1)}.
\end{equation}

We now determine $b_\lambda$ so that~\eqref{eq:linear comb} holds.
Notice that~\eqref{eq:goal}, for $1 \leq j \leq \Lambda$, can be written as the matrix equation
\begin{equation}
    \label{eq:matrix eqn}
    DS
    \begin{pmatrix}
         b_1\\
         b_2\\
         \vdots\\
         b_\Lambda
    \end{pmatrix}
    =
    \frac{(-1)^{\Lambda-1}}{(\Lambda-1)!}
    \begin{pmatrix}
        \,\,\,\, {\Lambda-1 \choose 0}\\
        -{\Lambda-1 \choose 1}\\
         \vdots\\
         (-1)^{\Lambda-1}{\Lambda-1 \choose \Lambda-1}
    \end{pmatrix},
\end{equation}
where $S$ is the upper triangular matrix
\begin{equation}
    \label{eq:S}
    S = \left( (-1)^{i+j}S(j,i) \right)_{\Lambda \times \Lambda},
\end{equation}
and $D$ is the diagonal matrix
\begin{equation}
    \label{eq:D}
    D=
    \begin{pmatrix}
        1!&  &  &&  \\
         & 2! & &&  \\
        &  & 3!&& \\
        &  & & \ddots&\\
         & &  & &\Lambda!
    \end{pmatrix}.
\end{equation}
Thus,
\begin{equation}
    \label{eq:matrix eqn b}
    \begin{pmatrix}
         b_1\\
         b_2\\
         \vdots\\
         b_\Lambda
    \end{pmatrix}
    =
    \frac{(-1)^{\Lambda-1}}{(\Lambda-1)!}
    S^{-1} D^{-1}
    \begin{pmatrix}
        \,\,\,\, {\Lambda-1 \choose 0}\\
        -{\Lambda-1 \choose 1}\\
         \vdots\\
         (-1)^{\Lambda-1}{\Lambda-1 \choose \Lambda-1}
    \end{pmatrix},
\end{equation}
with
\begin{equation}
    \label{eq:D inverse}
    D^{-1}=
    \begin{pmatrix}
        1!&  &  &&  \\
         & 1/2! & &&  \\
        &  & 1/3!&& \\
        &  & & \ddots&\\
         & &  & &1/\Lambda!
    \end{pmatrix}
\end{equation}
and
\begin{equation}
    \label{eq:S inverse}
    S^{-1}= \left( (-1)^{i+j}s(j,i) \right)_{\Lambda \times \Lambda},
\end{equation}
where $s(j,i)$ are the Stirling numbers of the first kind. Therefore, the solution
to equation~\eqref{eq:goal} is given by:
\begin{equation}
    \label{eq:b_lambda}
    b_\lambda = \frac{(-1)^{\Lambda+\lambda}}{(\Lambda-1)!}
    \sum_{j=\lambda}^{\Lambda} \frac{s(j,\lambda)}{j!}{\Lambda-1\choose \lambda-1},
\end{equation}
thus yielding~\eqref{eq:desired formula}.

\section{Expansion for Dirichlet $L$-functions}

Let $\chi$ be a non-trivial Dirichlet character for the modulus $q$.
Let $L(s,\chi)$ be the Dirichlet $L$-function
\begin{equation}
    \label{eq:L}
    L(s,\chi) := \sum_{n=1}^\infty \chi(n) n^{-s}.
\end{equation}
The defining sum is absolutely convergent in $\Re{s}>1$ and,
if $\chi$ is not the trivial character mod $q$, conditionally convergent
for $\Re{s}>0$. Rearranging the sum, for $\Re{s}>1$, we can write it as
\begin{equation}
    \label{eq:L in terms of hurwitz}
    \frac{1}{q^s} \sum_{m=1}^{q-1} \chi(m) \zeta(s,m/q).
\end{equation}
Applying~\eqref{eq:hurwitz expansion} and changing order of summation gives
\begin{equation}
    \label{eq:L expansion}
    L(s,\chi)
    =
    \frac{1}{\Gamma(s)}
    \frac{1}{q^s}
    \sum_{k=0}^\infty
    \alpha_k(s) \Gamma(s+k-1)
    \sum_{m=1}^{q-1}
    \frac{\chi(m) \Gamma(m/q)}{\Gamma(s+k+m/q-1)}.
\end{equation}
Similarly, for any positive integer $q$,
\begin{equation}
    \label{eq:zeta expansion}
    \Gamma(s) \zeta(s)
    = \frac{1}{q^s}
    \sum_{k=0}^\infty
    \alpha_k(s) \Gamma(s+k-1)
    \sum_{m=1}^{q}
    \frac{\Gamma(m/q)}{\Gamma(s+k+m/q-1)}.
\end{equation}
Notice that the sum over $m$ in~\eqref{eq:zeta expansion} runs up to $m\leq q$
rather than $q-1$ since, without the factor $\chi(m)$, we cannot ignore the
$m=q$ term.

Next, if we substitute $s=1$ into~\eqref{eq:L expansion} then only the $k=0$
term contributes because $\alpha_k(1)=0$ if $k\geq 1$, and
$\Gamma(s+k-1)/\Gamma(s+k+m/q-1)$ does not have a pole at $s=1$
when $k \geq 1$.

To deal with the $k=0$ term, consider, for
$a>0$, the Laurent series about $s=1$
\begin{equation}
    \frac{\Gamma(s-1)}{\Gamma(s+a-1)} = \frac{1}{\Gamma(a)(s-1)}
    - \frac{\gamma+\psi(a)}{\Gamma(a)}
    +O(s-1),
\end{equation}
where
\begin{equation}
    \psi(a) := \frac{\Gamma'(a)}{\Gamma(a)}.
\end{equation}
Substituting this, with $a=m/q$, into the $k=0$ term of~\eqref{eq:L expansion},
summing over $m$ and using $\sum \chi(m) = 0$, so that the terms involving $1/(s-1)$
and $\gamma$ sum to zero, and then letting $s=1$, we get the well known (see, for example,
Proposition 10.2.5 of~\cite{Co}) formula
\begin{equation}
    L(1,\chi) = \frac{-1}{q} \sum_{m=1}^{q-1} \chi(m) \frac{\Gamma'(m/q)}{\Gamma(m/q)}.
\end{equation}
Next, let $r$ be a non-negative integer. Consider the value of~\eqref{eq:L expansion}
at $s=-r$ which involves, on the rhs, the terms with $k=0,\ldots,r+1$.
Using the fact that the residue at $s=-r$ of $\Gamma(s+k-1)$ is equal to
$(-1)^{r-k+1}/(r-k+1)!$, we have that
\begin{equation}
    L(-r,\chi) = r! q^r
    \sum_{k=0}^{r+1} \frac{(-1)^{k-1}\alpha_k(-r)}{(r+1-k)!}
    \sum_{m=1}^{q-1}
    \frac{\chi(m) \Gamma(m/q)}{\Gamma(m/q+k-r-1)}.
\end{equation}
We can simplify slightly using $\Gamma(m/q)/\Gamma(m/q-r+k-1) = (m/q-1)\ldots(m/q-r+k-1)$
for $k<r+1$ and $1$ for $k=r+1$. The latter case does not contribute when summed against $\chi(m)$,
and we get
\begin{equation}
    \label{eq:L -r}
    L(-r,\chi) = r! q^r
    \sum_{k=0}^{r} \frac{(-1)^{k-1}\alpha_k(-r)}{(r+1-k)!}
    \sum_{m=1}^{q-1} \chi(m) (m/q-1)\ldots(m/q+k-r-1).
\end{equation}
For example, when $r=0$ the above gives
\begin{equation}
    L(0,\chi) =
    - \sum_{m=1}^{q-1} \chi(m) (m/q-1)
    = -\frac{1}{q} \sum_{m=1}^{q-1} \chi(m) m.
\end{equation}
The coefficients $\alpha_k(-r)$ that are needed for a given $r$ in~\eqref{eq:L -r}
are easily computed using, for instance, the recurrence~\eqref{eq:recur}.

We can write a similar formula for $\zeta(-r)$, though, without the $\chi(m)$,
we cannot ignore the $k=r+1$ term, and, furthermore, the sum over $m$ needs to be up to $m
\leq q$:
\begin{equation}
    \zeta(-r) = r! q^r
    \sum_{k=0}^{r+1} \frac{(-1)^{k-1}\alpha_k(-r)}{(r+1-k)!}
    \sum_{m=1}^{q} (m/q-1)\ldots(m/q+k-r-1).
\end{equation}
The innermost summand above is taken to equal 1 if $k=r+1$. For example, if $q=1$,
we get
\begin{equation}
    \zeta(-r) = (-1)^r r! \alpha_{r+1}(-r).
\end{equation}
In Section 7 of~\cite{R}, the author showed that $\alpha_{r+1}(-r) =
B_{r+1}/(r+1)!$. Alternatively, one can use the $\lambda=1$ case of~\eqref{eq:3 intro},
and the association~\eqref{eq:alpha vs B} to compute the residue at $s=-r+1$ from the $k=r+1$ term
of~\eqref{eq:3 intro}, with special attention paid to the case $r=0$.
Thus, the above gives Euler's formula $\zeta(-r) = (-1)^r B_{r+1}/(r+1)$.

Next, we develop a second formula for $L(s,\chi)$ evaluated at negative integers.
Using~\eqref{eq:hurwitz zeta with lambda shift} we have, for positive integer $\lambda$:
\begin{eqnarray}
    \label{eq:L with shift}
    L(s-\lambda,\chi)
    &=&
    \frac{1}{\Gamma(s)}
    \frac{1}{q^{s-\lambda}}
    \sum_{k=0}^\infty
    \alpha_k(s)
    \sum_{m=1}^{q-1}
    \chi(m) \Gamma(m/q) \notag \\
    &\times& \sum_{j=0}^\lambda
    \frac{\Gamma(s+k-j-1)}{\Gamma(s+k+m/q-j-1)} c_{m/q}(\lambda,j). \notag \\
\end{eqnarray}
Substituting $s=1$, we get contributions from the terms $k=0$, and from $j=k$ with
$j=1,\ldots, \lambda$. The latter is accounted for by the pole of $\Gamma(s+k-j-1)$
balancing out against the zero of $\alpha_k(s)$ at $s=1$. Thus
\begin{eqnarray}
    \label{eq:L formula 2}
    L(1-\lambda,\chi)
    &=& q^{\lambda-1}
    \sum_{j=0}^\lambda
    \left(
        \frac{(-1)^j}{j!}
        + \sum_{k=1}^j
        \frac{(-1)^{k-j}}{(j-k)!} \alpha_k'(1)
    \right) \notag \\
    &\times& \sum_{m=1}^{q-1}
    \chi(m) c_{m/q}(\lambda,j) (m/q-1)\ldots(m/q-j+k). \notag \\
\end{eqnarray}

\section{General Dirichlet series}

Let
\begin{equation}
    L(s) = \sum_{n=1}^\infty \frac{b(n)}{n^s}
    \label{eq:dirichlet series}
\end{equation}
be a Dirichlet series, assumed to be absolutely convergent for $\Re{s}>\sigma_1$.

Substitute formula~\eqref{eq:gamma}, with $w=n$
into the Dirichlet series, and change order of summation over $n$ and $k$:
\begin{eqnarray}
    \label{eq:hurwitz L expansion}
    L(s)
    &=&
    \frac{1}{\Gamma(s)}
        \sum_{n=1}^\infty b(n)
        \sum_{k=0}^\infty
        \frac{\alpha_k(s) \Gamma(s+k)(n-1)!}{\Gamma(s+k+n)} \notag \\
    &=&
    \frac{1}{\Gamma(s)}
        \sum_{k=0}^\infty \alpha_k(s)
        \sum_{n=1}^\infty b(n) \frac{\Gamma(s+k)(n-1)!}{\Gamma(s+k+n)}.
\end{eqnarray}
Changing order of the sums can be justified for $\Re{s}> \sigma_1$, but we do not
do so here, because our emphasis below will mainly be on the first equality, rather than the second.
However, before turning our attention to the first equality, we
note that, if $b(n)$ is identically equal to 1 for $n\geq N+1$ and 0 otherwise, then the
inner sum over $n$ in the last line above equals
\begin{equation}
    \frac{\Gamma(s+k-1)N!}{\Gamma(s+k+N)}.
\end{equation}
This can be seen by taking the identity
\begin{eqnarray}
    \label{eq:hurwitz L manipulate}
    \int_0^1 t^{s+k-1} \sum_{n=N+1}^\infty (1-t)^{n-1} dt
    = \sum_{n=N+1}^\infty \frac{\Gamma(s+k) \Gamma(n)}{\Gamma(s+k+n)},
\end{eqnarray}
and summing the geometric series on the lhs. Thus, in this case,
equation~\eqref{eq:hurwitz L expansion} reduces to~\eqref{eq:zeta via hurwitz b}. The fact that we can
simplify the sum over $n$ is what makes~\eqref{eq:hurwitz L expansion}
particularly useful in the
case of the Riemann and Hurwitz zeta functions, and Dirichlet $L$-functions.

Next we explore the first equality of~\eqref{eq:hurwitz L expansion}.
Using the partial
fraction expansion~\eqref{eq:parfrac}, we have
\begin{equation}
    \label{eq:parfrac b}
    \frac{(n-1)!}{\Gamma(s+k+n)}
    =\sum_{m=0}^{n-1} \frac{(-1)^m}{s+k+m} {n-1 \choose m}.
\end{equation}
Therefore,
\begin{equation}
    \label{eq:L manipulate a}
    L(s) = \frac{1}{\Gamma(s)}
    \sum_{n=1}^\infty b(n)
    \sum_{k=0}^\infty \alpha_k(s)
    \sum_{m=0}^{n-1} \frac{(-1)^m}{s+k+m} {n-1 \choose m}.
\end{equation}
But bound~\eqref{eq:alpha bound} allows us to rearrange the two inner sums over $k$ and $m$:
\begin{equation}
    \sum_{k=0}^\infty \alpha_k(s) \sum_{m=0}^{n-1} \frac{(-1)^m}{s+k+m} {n-1 \choose m}
    =
    \sum_{m=0}^{n-1} (-1)^m {n-1 \choose m}
    \sum_{k=0}^\infty \frac{\alpha_k(s)}{s+k+m}.
\end{equation}
The resulting inner sum over $k$ equals
\begin{eqnarray}
    \label{eq:binomial via int}
    \int_0^1 (-\log(1-t)/t)^{s-1} t^{s+m-1} dt
    &=&
    \int_0^\infty x^{s-1} \exp(-x) (1-\exp(-x))^m dx \notag \\
    &=&
    \Gamma(s) \sum_{j=0}^m \frac{(-1)^j}{(j+1)^s} {m \choose j},
\end{eqnarray}
which we see using~\eqref{eq:1}, then changing variable $t=1-\exp(-x)$,
expanding $(1-\exp(-x))^m$ using the binomial theorem,
and integrating termwise.
We thus have the identity
\begin{equation}
    \label{eq:alpha binomial identity}
    \frac{1}{\Gamma(s)}
    \sum_{k=0}^\infty \frac{\alpha_k(s)}{s+k+m}
    =
    \sum_{j=0}^m \frac{(-1)^j}{(j+1)^s} {m \choose j}.
\end{equation}
Notice that, while in the integrals in~\eqref{eq:binomial via int}
we require $\Re{s}>-m$, both sides of the above equation are entire functions of $s$,
and hence, by analytic continuation, equality holds for all $s$.

Hence,
\begin{equation}
    \label{eq:L manipulate b}
    L(s) =
    \sum_{n=1}^\infty b(n)
    \sum_{m=0}^{n-1} (-1)^m {n-1 \choose m}
    \sum_{j=0}^m \frac{(-1)^j {m \choose j}}{(j+1)^s}.
\end{equation}
Since the inner double sum is just another way
to express $1/n^s$, formula~\eqref{eq:L manipulate b} holds for all $s$ for which the
defining Dirichlet series converges conditionally.

We can also arrive at the above more directly by starting with the binomial coefficient
identity
\begin{equation}
    \label{eq:binomial identity}
    \sum_{m=j}^{n-1} (-1)^m {n-1 \choose m} {m \choose j} =
    \begin{cases}
        (-1)^{n-1}, &\text{if $j=n-1$,} \\
        0, &\text{otherwise,}
    \end{cases}
\end{equation}
for $j\geq 0$, and $n\geq 1$. To prove this identity, consider the coefficient of $z^{n-1}$
in $(z-1)^{n-1} \sum_{m=j}^\infty {m \choose j} z^m = (-1)^{n-1} z^j (1-z)^{n-j-2}$.

Thus
\begin{equation}
    \label{eq:n^s identity}
    \sum_{j=0}^{n-1} \frac{(-1)^j}{(j+1)^s}
    \sum_{m=j}^{n-1} (-1)^m {n-1 \choose m} {m \choose j} =
    \frac{1}{n^s},
\end{equation}
because only the $j=n-1$ terms survives with an extra $(-1)^{n-1}$ introduced from the inner sum.
Substituting into the Dirichlet series for $L(s)$, and rearranging the inner double sum
gives~\eqref{eq:L manipulate b}.

\subsection{Heuristic manipulations}

Next we proceed heuristically, but we will justify our resulting formulas in certain cases
in the next section.

Changing order of summation, and ignoring for the moment the fact that the
resulting sums over $n$ diverge unless $b(n)$ is rapidly decreasing, gives
\begin{equation}
    \label{eq:L rearrange n}
    L(s) =
    \sum_{m=0}^{n-1} \frac{(-1)^m}{m!}
    \sum_{n=1}^\infty b(n) (n-1) \ldots (n-m)
    \sum_{j=0}^m \frac{(-1)^j {m \choose j}}{(j+1)^s}.
\end{equation}
We interpret the typically divergent sum over $n$ as being equal to the $m$-th derivative,
evaluated at $z=1$, of the analytic continuation, if it exists, of the series
\begin{equation}
    \label{eq:f}
    f(z) := \sum_{n=1}^\infty b(n) z^{n-1},
\end{equation}
thus suggesting
\begin{equation}
    \label{eq:L rearrange n b}
    L(s) =
    \sum_{m=0}^{n-1} \frac{(-1)^m f^{(m)}(1)}{m!}
    \sum_{j=0}^m \frac{(-1)^j {m \choose j}}{(j+1)^s}.
\end{equation}

For example, let $L(s)=\zeta(s)(1-2/2^s) = 1-2^{-s}+3^{-3}-4^{-s}+\ldots$ be the
alternating zeta function. We consider here the alternating zeta function, rather than
$\zeta$ itself, because
the pole of $\zeta(s)=1$ creates a slight complication that we wish to
circumvent.
The factor $(1-2/2^s)$ has a zero at $s=1$, balancing the
pole of $\zeta(s)$. Thus, $\zeta(s)(1-2/2^s)$ extends to an entire function.
Then $f(z)=1-z+z^2-z^3+\ldots= 1/(1+z)$, and $f^{(m)}(1)= (-1)^m m!/2^{m+1}$,
hence~\eqref{eq:L rearrange n} reduces to
\begin{equation}
    \label{eq:zeta alt}
    \zeta(s)(1-2/2^s) =
    \sum_{m=0}^\infty \frac{1}{2^{m+1}}
    \sum_{j=0}^m \frac{(-1)^j {m \choose j}}{(j+1)^s},
\end{equation}
i.e. to Hasse's formula described in the introduction.

We can also, in~\eqref{eq:L rearrange n}, expand
\begin{equation}
    \label{eq:stirling expand}
    (n-1)\ldots(n-m)
    =\frac{1}{n} n(n-1)\ldots(n-m)
    =
    \frac{1}{n} \sum_{l=0}^{m+1} s(m+1,l) n^l.
\end{equation}
Here, we have multiplied and divided, in the first equality, by $n$ so as to
have, in the second equality, a polynomial in $n$ rather than in $n-1$.
Substitute the rhs of~\eqref{eq:stirling expand} into~\eqref{eq:L rearrange n},
and rearrange the sum over $n$ and $l$
to get
\begin{equation}
    \label{eq:L manipulate d}
    L(s) =
    \sum_{m=0}^\infty \frac{(-1)^m}{m!}
    \sum_{l=0}^{m+1} s(m+1,l)
    \sum_{n=1}^\infty \frac{b(n)}{n^{1-l}}
    \sum_{j=0}^m \frac{(-1)^j {m \choose j}}{(j+1)^s}.
\end{equation}
We interpret the typically divergent sum $\sum_{n=1}^\infty b(n)n^{l-1}$ as being equal to
$L(1-l)$, suggesting the following interpolation formula for $L(s)$:
\begin{equation}
    \label{eq:L interpolate}
    L(s) =
    \sum_{m=0}^\infty \frac{(-1)^m}{m!}
    \sum_{l=0}^{m+1} s(m+1,l) L(1-l)
    \sum_{j=0}^m \frac{(-1)^j {m \choose j}}{(j+1)^s}.
\end{equation}
Now, $s(m+1,0)=0$ for all $m\geq 0$, and thus we can start the sum
at $l=1$ unless $L(s)$ has a pole at $s=1$.

Consider again the case $L(s)=\zeta(s)(1-2/2^s)$.
In this case, for $l\geq 1$,
\begin{equation}
    \label{eq:eta knopp}
    L(1-l)= \zeta(1-l)(1-2^l) = (-1)^{l-1} (1-2^l) B_l/l,
\end{equation}
and the sum over $m$ in~\eqref{eq:L interpolate} equals
\begin{equation}
    \sum_{l=1}^{m+1} (-1)^{l-1} (1-2^l) s(m+1,l) B_l/l.
    \label{eq:funky 1}
\end{equation}
Comparing with~\eqref{eq:zeta alt} suggests that the above equals
\begin{equation}
    (-1)^m m!/2^{m+1}.
    \label{eq:funky 2}
\end{equation}
In the next section we will rigorously prove the above manipulations in a few
cases, including the alternating zeta function, and thus yielding equality
between~\eqref{eq:funky 1} and~\eqref{eq:funky 2}.

Note that we can also consider the above manipulations for the Dirichlet series with coefficients
$b(n)/n^{s_0}$ for some $s_0 \in \C$. Writing
\begin{equation}
    \label{eq:L with s_0}
    L(s+s_0) = \sum_{n=1}^\infty \frac{b(n)}{n^{s_0}} \frac{1}{n^{s}}
\end{equation}
and letting
\begin{equation}
    \label{eq:L s_0}
    L(s_0,z) := \sum_{n=1}^\infty \frac{b(n)}{n^{s_0}} z^{n-1},
\end{equation}
we anticipate (in certain cases) two formulas:
\begin{equation}
    \label{eq:L s_0 a}
    L(s+s_0) =
    \sum_{m=0}^{n-1} \frac{(-1)^m L^{(m)}(s_0,1)}{m!}
    \sum_{j=0}^m \frac{(-1)^j {m \choose j}}{(j+1)^s},
\end{equation}
where $L^{(m)}(s_0,1)$ is the $m$-th derivative with respect to $z$ of $L(s_0,z)$ evaluated
at $z=1$, and
\begin{equation}
    \label{eq:L s_0 b}
    L(s+s_0) =
    \sum_{m=0}^\infty \frac{(-1)^m}{m!}
    \sum_{l=0}^{m+1} s(m+1,l) L(s_0+1-l)
    \sum_{j=0}^m \frac{(-1)^j {m \choose j}}{(j+1)^s}.
\end{equation}

\subsection{Justifying the manipulations}

We are able to prove the formulas of the previous section, specifically
equations~\eqref{eq:L s_0 a} and~\eqref{eq:L s_0 b} for the alternating zeta
and alternating Hurwitz zeta functions, and Dirichlet $L$-functions $L(s,\chi)$,
where $\chi$ is a non-trivial character for the modulus $q\leq 5$. In
Section~\ref{sec:summation dirichlet} we develop a related summation formula
that can be used to give interesting formulas valid for all $q$.

For higher degree $L$-functions, for example, degree 2 $L$-functions associated
to cusp forms of given weight, level, and character, it appears from numerical
experiments that formula~\eqref{eq:L s_0 b} holds but with some adjustment-
additional convergence producing terms, arising from a smoothed approximate
functional equation, seem to be needed. We will revisit the issue of higher
degree $L$-functions in a future paper.

We start by modifying equation~\eqref{eq:L manipulate b} by introducing an
extra parameter $z$. Let
\begin{equation}
    \label{eq:L z}
    L(s,z) :=
    \sum_{n=1}^\infty \frac{b(n)}{n^s} z^{n-1}.
\end{equation}
Assume that there exists a real number $r$ such that
\begin{equation}
    \label{eq:b}
    b(n) \ll n^{r},
\end{equation}
so that the coefficients $b(n)$ do not grow too quickly, and the series
defining $L(s,z)$ converges for all $s \in \C$ and $|z|<1$.
Then, as before, replacing $1/n^s$ by~\eqref{eq:n^s identity}, we rigorously have
\begin{equation}
    \label{eq:L z b}
    L(s,z) =
    \sum_{n=1}^\infty b(n) z^{n-1}
    \sum_{m=0}^{n-1} (-1)^m {n-1 \choose m}
    \sum_{j=0}^m \frac{(-1)^j {m \choose j}}{(j+1)^s}.
\end{equation}
Now, for $|z|< 1/2$, the sum over $n$ and $m$ converges absolutely. That is
because~\eqref{eq:j sum asympt} shows that the sum over $j$ is bounded, in the $m$-aspect,
by an amount which is generously $O(1)$, with the implied
constant uniform for $s$ in compact subsets of $\C$. Therefore, the absolute value
of the rhs of the above is, on using $\sum_{m=0}^{n-1} {n-1 \choose m} = 2^{n-1}$,
\begin{equation}
    \ll_s
    \sum_{n=1}^\infty |b(n)| |2z|^{n-1}.
    \label{eq:L s z abs}
\end{equation}
This converges for all $|z|<1/2$ because $b(n)$ is assumed to grow at most
polynomially in $n$. Therefore, for $|z|<1/2$, and $s \in \C$, we have on
rearranging the sums over $m$ and $n$ and pulling out a $z^m/m!$ from the sum
over $n$:
\begin{equation}
    \label{eq:L z c}
    L(s,z) =
    \sum_{m=0}^\infty \frac{(-z)^m}{m!}
    \sum_{n=m+1}^\infty b(n) (n-1)\ldots(n-m) z^{n-m-1}
    \sum_{j=0}^m \frac{(-1)^j {m \choose j}}{(j+1)^s}.
\end{equation}
Our next goal is to obtain the analytic continuation in $z$ of the sum over $n$
and to substitute $z=1$. In general, this cannot be done, but in some
examples we are able to do so.

\subsection{Alternating zeta function}
\label{sec:alt zeta}

The alternating zeta function is defined to be
\begin{equation}
    \label{eq:alt zeta}
    \eta(s)
    := \sum_{n=1}^{\infty} \frac{(-1)^{n-1}}{n^{s}}.
\end{equation}
We wish to allow shifts, so we let $b(n)=(-1)^{n-1}/n^{s_0}$ for some $s_0 \in \C$.
Thus, we are considering
\begin{equation}
    \label{eq:alt zeta s z}
    \eta(s+s_0,z) := \sum_{n=1}^{\infty} \frac{(-1)^{n-1}}{n^{s_0}} \frac{z^{n-1}}{n^s}.
\end{equation}
This can also be regarded as a special case of the Lerch zeta function.
We have introduced the parameter $s_0$ since we wish to examine the general
expansions of the form considered in~\eqref{eq:L s_0 b}.

In this case, the sum over $n$ in~\eqref{eq:L z c} reduces to
\begin{equation}
    \label{eq:alt zeta diff}
    \sum_{n=m+1}^\infty \frac{(-1)^{n-1}}{n^{s_0}} (n-1)\ldots(n-m) z^{n-m-1}.
\end{equation}
Now
\begin{equation}
    \sum_{n=1}^{\infty} \frac{(-1)^{n-1}}{n^{s_0}} z^{n-1}
    =
    \frac{1}{\Gamma(s_0)} \int_0^\infty \frac{x^{s_0-1}\exp(-x)}{1+z\exp(-x)} dx,
    \label{eq:integral eta}
\end{equation}
which we can prove by expanding the denominator as a geometric series and
integrating termwise. The above integral converges and is analytic in both variables
for $\Re{s_0} > 0$ and
$z \notin (-\infty,-1]$. Differentiating the lhs termwise $m$ times
within its disc of convergence
gives~\eqref{eq:alt zeta diff}. Hence differentiating the rhs $m$ times yields
\begin{eqnarray}
    \label{eq:alt zeta diff b}
    &&\sum_{n=m+1}^\infty \frac{(-1)^{n-1}}{n^{s_0}} (n-1)\ldots(n-m) z^{n-m-1} \notag \\
    &=&  \frac{(-1)^m m!}{\Gamma(s_0)}
    \int_0^\infty \frac{x^{s_0-1}\exp(-x(m+1))}{(1+z\exp(-x))^{m+1}} dx.
\end{eqnarray}
Also note that we can bound these derivatives using knowledge of the radius of convergence
of the Taylor series, about the point $w=z$ of the function $\eta(s_0,w)$. Consider the Taylor
expansion
\begin{equation}
    \eta(s_0,w) = \sum_{m=0}^\infty \frac{\eta^{(m)}(s_0,z)}{m!}(w-z)^m.
    \label{eq:taylor}
\end{equation}
Here $\eta^{(m)}$ refers to the $m$-th derivative of
$\eta$ with respect to the second variable.
Assume that $\Re{z}>-1$, so that the nearest singularity of the rhs
of~\eqref{eq:integral eta} is the point $-1$. Thus, the above converges
absolutely for all $|w-z|<|1+z|$, in particular, at $w=0$ if $|z|<|1+z|$. The
latter holds when $\Re{z}>-1/2$ (as can be seen by drawing $z$ and $1+z$).

Thus, returning to~\eqref{eq:L z c} (with $L(s)=\eta(s)$),
for $\Re{z}>-1/2$, $\Re{s_0}>0$, $s \in \C$,
\begin{equation}
    \eta(s+s_0,z)=
    \frac{1}{\Gamma(s_0)}
    \sum_{m=0}^\infty z^m
    \left(
        \int_0^\infty \frac{x^{s_0-1}\exp(-x(m+1))}{(1+z\exp(-x))^{m+1}} dx
    \right)
    \sum_{j=0}^m \frac{(-1)^j {m \choose j}}{(j+1)^s}.
    \label{eq:alt zeta expansion}
\end{equation}
(recall, for the purpose of convergence, that the sum over $j$ is $O_s(1)$).
Therefore, substituting $z=1$ we get
\begin{equation}
    \label{eq:alt zeta analytic continuation}
    \eta(s+s_0) =
    \frac{1}{\Gamma(s_0)}
    \sum_{m=0}^\infty
    \left(
        \int_0^\infty \frac{x^{s_0-1}\exp(-x(m+1))}{(1+\exp(-x))^{m+1}} dx
    \right)
    \sum_{j=0}^m \frac{(-1)^j {m \choose j}}{(j+1)^s},
\end{equation}
thus justifying~\eqref{eq:L s_0 a} in the case of $L(s)=\eta(s)$.
For any given $s_0$ satisfying $\Re{s_0}>0$, the rhs above converges
uniformly in $s$ on compact subsets of $\C$, thus also giving the analytic continuation
of the alternating zeta function.

We can extend the validity of~\eqref{eq:alt zeta analytic continuation} in the $s_0$ aspect
by integrating by parts
\begin{equation}
    \int_0^\infty \frac{x^{s_0-1}\exp(-x(m+1))}{(1+\exp(-x))^{m+1}} dx
    = \frac{(m+1)}{s_0}
    \int_0^\infty \frac{x^{s_0}\exp(-x(m+1))}{(1+\exp(-x))^{m+2}} dx.
    \label{eq:int by parts}
\end{equation}
While the integral on the lhs is convergent for $\Re{s_0}>0$, the integral on
the right converges for $\Re{s_0}>-1$. Therefore, for $\Re{s_0}>-1$ and $s \in \C$,
\begin{eqnarray}
    \label{eq:alt zeta analytic continuation b}
    &&\eta(s+s_0) = \notag \\
    &&\frac{1}{\Gamma(s_0+1)}
    \sum_{m=0}^\infty (m+1)
    \left(
        \int_0^\infty \frac{x^{s_0}\exp(-x(m+1))}{(1+\exp(-x))^{m+2}} dx
    \right)
    \sum_{j=0}^m \frac{(-1)^j {m \choose j}}{(j+1)^s}.\notag \\
\end{eqnarray}
For example, substituting $s_0=0$ gives Hasse's formula~\eqref{eq:hasse} as
described in the introduction.

Because the $m$-th term decreases exponentially in $m$ for given $s_0$ (with
$\Re{s_0}>-1$), and uniformly for $s$ in compacts subsets of $\C$,
equation~\eqref{eq:alt zeta analytic continuation b} also gives the analytic
continuation of $\eta$.

The above also shows that, for {\it any} $s \in \C$ and $\Re{s_0}>-1$, the limit, as $z \to 1$ of
the rhs of~\eqref{eq:alt zeta s z} exists and is equal to $\eta(s+s_0)$ as expressed,
for example, in~\eqref{eq:alt zeta analytic continuation b}. Therefore, returning
to~\eqref{eq:alt zeta diff} and expanding using~\eqref{eq:stirling expand} we get
\begin{eqnarray}
    \label{eq:alt zeta analytic continuation c}
    &&\eta(s+s_0) = \notag \\
    &&
    \sum_{m=0}^\infty \frac{(-1)^m}{m!}
    \sum_{l=1}^{m+1} s(m+1,l) \eta(s_0+1-l)
    \sum_{j=0}^m \frac{(-1)^j {m \choose j}}{(j+1)^s}. \notag \\
\end{eqnarray}
While we initially assumed $\Re{s_0}>-1$, we can extend~\eqref{eq:alt zeta
analytic continuation b} to $\Re{s_0}>-M$, where $M$ is any positive integer,
by repeatedly integrating by parts. Therefore, equation~\eqref{eq:alt zeta
analytic continuation c} holds for all $s_0,s \in \C$.

\subsection{Dirichlet $L$-functions}
\label{sec:dirichlet hasse}

Let $\chi$ be a non-trivial Dirichlet character for the modulus $q$. Define
\begin{equation}
    \label{eq:L s z chi}
    L(s_0,z,\chi) = \sum_{n=1}^\infty \frac{\chi(n)}{n^{s_0}} z^{n-1}, \qquad |z|<1.
\end{equation}
Now, for $\Re{s_0}>0$,
\begin{equation}
    \label{eq:L s z chi b}
    L(s_0,z,\chi) =
    \frac{1}{\Gamma(s_0)} \int_0^\infty \frac{x^{s_0-1}\sum_{a=1}^{q-1} \chi(a) z^{a-1}\exp(-ax)}{1-z^q\exp(-qx)} dx,
\end{equation}
which we can see by expanding the denominator as a geometric series and using
the periodicity mod $q$ of $\chi$. Note that, for given $x$, the singularity of
the integrand at $z=\exp(x)$ is removable because the numerator vanishes at
$z=\exp(x)$. Therefore, the rhs is analytic in some neighbourhood, depending on
$q$, of $z=1$.

As in the discussion concerning $\eta$, we can get an estimate for the size of
the $m$-th derivative with respect to $z$ of $L(s_0,z,\chi)$, evaluated at
$z=1$ by considering the singularities of the rhs of~\eqref{eq:L s z chi b}
which occur along rays emanating outward from the non-one $q$-th roots of unity, i.e.
at $z=r\exp(2\pi i j/q)$, $r>1$, $1 \leq j < q$.

We require that the distance
from the point 1 to the nearest singularity, $\exp(\pm 2\pi i/q)$ be greater
than 1 so that the sum analogous to~\eqref{eq:alt zeta expansion} continues analytically to $z=1$.
This occurs when $q<5$ (at $q=6$, $\exp(2\pi i/6)$ has a distance of 1
from the point 1, and for $q>6$ the closest root is even closer). Thus, for any
non-trivial $\chi$ mod $q\leq 5$ we have:
\begin{equation}
    \label{eq:L s chi interpolate}
    L(s+s_0,\chi) =
    \sum_{m=0}^\infty \frac{(-1)^m}{m!}
    \sum_{l=1}^{m+1} s(m+1,l) L(s_0+1-l)
    \sum_{j=0}^m \frac{(-1)^j {m \choose j}}{(j+1)^s},
\end{equation}
valid for all $s,s_0 \in \C$.

\section{A summation formula}
\label{sec:summation}

Next we consider a generalization of~\eqref{eq:n^s identity}. The same binomial
coefficient identity which gave that formula for $1/n^s$ can be used for any function
$h(n)$ on the positive integers:
\begin{equation}
    \label{eq:h identity}
    h(n)=
    \sum_{j=0}^{n-1} (-1)^{j} h(j+1)
    \sum_{m=j}^{n-1} (-1)^m {n-1 \choose m} {m \choose j},
\end{equation}
because only the $j=n-1$ term survives. Therefore rearranging sums,
\begin{equation}
    \label{eq:h identity b}
    h(n)=
    \sum_{m=0}^{n-1} (-1)^m {n-1 \choose m}
    \sum_{j=0}^{m} (-1)^{j} h(j+1) {m \choose j}.
\end{equation}

Let $g,h: \N \to \C$ be two functions on the positive integers.
Apply the above identity for $h(n)$ to the power series
\begin{equation}
    \label{eq:gh}
    \sum_{n=1}^\infty g(n)h(n) z^{n-1},
\end{equation}
and then rearrange the sum over $m$ and $n$, pulling out a $1/m!$ and $z^m$
from the sum over $m$:
\begin{eqnarray}
    \label{eq:gh summation}
    &&\sum_{n=1}^\infty g(n)h(n) z^{n-1}
    \notag \\
    &&=\sum_{m=0}^\infty \frac{(-z)^m}{m!}
    \sum_{n=m+1}^\infty g(n) (n-1) \ldots (n-m) z^{n-m-1} \notag \\
    &&\qquad \qquad \qquad \qquad \times \sum_{j=0}^m (-1)^j h(j+1) {m \choose j}. \notag \\
\end{eqnarray}
The latter step of rearranging the sums over $m$ and $n$ can be justified for sufficiently small $z$
assuming modest bounds on $g,h$. For example, assume that $g(n)$ and $h(n)$ grow at most exponentially,
i.e. both are $O(\exp(cn))$, for some $c\geq 0$. Then, substituting~\eqref{eq:h identity b}
into the lhs of~\eqref{eq:gh summation}, the sum is, in absolute value,
\begin{equation}
    \ll \sum_{n=1}^\infty \exp(2cn) 3^{n} |z|^{n-1}.
\end{equation}
The $3^{n}$ results from applying the trivial bounds
$\sum_{j=0}^m  {m \choose j} = 2^m$, and $\sum_{m=0}^{n-1} {n-1
\choose m} 2^m = 3^{n-1}$ to the two sums in~\eqref{eq:h identity b},
while the $\exp(2cn)$ arises from the assumed bounds on $g$ and $h$.
Therefore, for $|z|< \exp(-2c)/3$ the above sum converges, and we may rearrange the
sums over $m$ and $n$.


We can view the sum over $n$ as being the $m$-th derivative with respect to $z$
of the power series
\begin{equation}
    \label{eq:G}
    G(z):=\sum_{n=1}^\infty g(n) z^{n-1}.
\end{equation}
Thus equation~\eqref{eq:gh summation} can be written as:
\begin{equation}
    \label{eq:gh summation b}
    \sum_{n=1}^\infty g(n)h(n) z^{n-1}
    =
    \sum_{m=0}^\infty G^{(m)}(z) \frac{(-z)^m}{m!}
    \sum_{j=0}^m (-1)^j h(j+1) {m \choose j}.
\end{equation}

Note that, in the previous sections, with the choice $h(n)=1/n^s$,
and the asymptotic formula~\eqref{eq:j sum asympt}, so that
$\sum_{j=0}^m (-1)^j {m \choose j}/(j+1)^s
\ll_s 1$, we were thus allowed to initially take $|z|<1/2$. The second
ingredient for $\eta$ and $L(s,\chi)$, with $g(n)=(-1)^{n-1}/n^{s_0}$ or
$\chi(n)/n^{s_0}$, was to obtain the analytic continuation of the rhs
of~\eqref{eq:gh summation} to a neighbourhood of the point $z=1$.
We achieved this by expressing the appropriate $G(z)$ in closed form
in~\eqref{eq:integral eta} and~\eqref{eq:L s z chi b}
from which we were able to deduce the rate of decay of $G^{(m)}(z)$.

The inner sum over $j$,
\begin{equation}
    \label{eq:finite diff}
    \sum_{j=0}^m (-1)^j h(j+1) {m \choose j},
\end{equation}
can be interpreted in terms of the finite difference operator. Define
\begin{eqnarray}
    \label{eq:delta}
    \Delta h(j) &=&  h(j+1)-h(j) \notag \\
    \Delta^m h(j) &=&  \Delta^{m-1} h(j+1) - \Delta^{m-1} h(j).
\end{eqnarray}
Then~\eqref{eq:finite diff} is equal to
\begin{equation}
    (-1)^m \Delta^m h(1),
\end{equation}
hence
\begin{equation}
    \label{eq:gh summation c}
    \sum_{n=1}^\infty g(n)h(n) z^{n-1}
    =
    \sum_{m=0}^\infty G^{(m)}(z) \Delta^{m} h(1) \frac{z^m}{m!}.
\end{equation}

Note that, because our application is to Dirichlet series, we prefer, above, to
write the $n$-the term of the lhs as $g(n)h(n) z^{n-1}$. However, in other applications,
as to the Hurwitz zeta function below,
one might prefer to write the $n$-th term as $g(n)h(n) z^{n}$ and start the sum at $n=0$.
Thus, as a variant, for functions $g,h$ on the non-negative integers
\begin{equation}
    \label{eq:gh summation d}
    \sum_{n=0}^\infty g(n)h(n) z^n
    =
    \sum_{m=0}^\infty \tilde{G}^{(m)}(z) \Delta^{m} h(0) \frac{z^m}{m!},
\end{equation}
where
\begin{equation}
    \label{eq:tilde G}
    \tilde{G}(z):=\sum_{n=0}^\infty g(n) z^{n}.
\end{equation}
The range of validity in $z$ of~\eqref{eq:gh summation c} and~\eqref{eq:gh summation d}
depends on the rates of decay
of $g(n)$, and $\Delta^{m} h(0)$ as discussed above, for example being valid on some disc centred at $0$
if $g(n)$ and $h(n)$ grow at most exponentially in $n$. The analytic continuation of the rhs to a given
point $z$ depends further on the location of the singularities of $G(w)$ or of $\tilde{G}(w)$
in relation to the point $z$.

To illustrate our summation formula, let $\lambda>0$, and take $g(n)=(-1)^{n-1} \lambda^n$,
$h(n) =n^{-s} \lambda^{-n}$. Then $G(z)= \lambda/(1+\lambda z)$, and
$G^{(m)}(1) (-1)^m/m! = \lambda^{m+1}/(1+\lambda)^{m+1}$. Formally, this choice of $g$ and $h$
gives Amore's generalization of Hasse's identity~\cite{A}:
\begin{eqnarray}
    \label{eq:amore}
    \eta(s) =
    \sum_{m=0}^\infty \frac{\lambda^{m+1}}{(1+\lambda)^{m+1}}
    \sum_{j=0}^m \frac{(-1)^j {m \choose j}}{(j+1)^s \lambda^{j+1}}. \notag \\
\end{eqnarray}
In order to justify the above, one also needs to obtain a rate of decay for the
sum over $j$ (depending on $\lambda$), and this is done in~\cite{A}.

Many variants are possible. For example, let $\lambda>0$, $g(n)=$ $(-1)^{n-1}$ $\lambda^{n-1}/(n-1)!$,
$h(n) =n^{-s} \lambda^{-n+1} (n-1)!$. Then $G(z)= \exp(-\lambda z)$, and
$G^{(m)}(1)=(-\lambda)^m \exp(-\lambda)$. Therefore, substituting into our summation formula,
and simplifying the factorials and binomial coefficient that appear, we get
\begin{equation}
    \label{eq:amore variant}
    \eta(s) = \exp(-\lambda) \sum_{m=0}^\infty \lambda^m
    \sum_{j=0}^m \frac{(-1)^j}{(j+1)^s \lambda^j (m-j)!}.
\end{equation}
We justify this formula for $\Re{s}>0$.
Note that, while $h(n)$ grows very quickly with $n$, to justify the above formula
what we need, for the purpose, here, of letting $z \to 1^-$
in~\eqref{eq:gh summation}, is an estimate for the size of the sum
\begin{equation}
    \label{eq:j sum amore variant}
    \sum_{j=0}^m \frac{(-1)^j}{(j+1)^s \lambda^j (m-j)!}.
\end{equation}
For $\Re{s}>0$, this sum is equal to
\begin{eqnarray}
    \label{eq:amore variant generating function}
    &&\frac{1}{\Gamma(s)} \int_0^\infty x^{s-1}
    \sum_{j=0}^m \frac{(-1)^j\exp(-(j+1)x)}{\lambda^j (m-j)!}.
\end{eqnarray}
Note that, as $m \to \infty$, the main contribution to the above integral
arises from smaller $x$, which we can see by breaking the integral into, say,
$\int_0^1+\int_1^\infty$, and on the second integral, breaking the sum into
$0<j\leq m/2$ and $m/2 <j \leq m$. We thus truncate the integral at $x=1$, and
then reverse the sum over $j$ to get, as $m \to \infty$,
\begin{eqnarray}
    \label{eq:amore variant generating function b}
    &&\sum_{j=0}^m \frac{(-1)^j}{(j+1)^s \lambda^j (m-j)!}
    \sim \notag \\
    &&\frac{(-1)^m}{\lambda^m\Gamma(s)} \int_0^1 x^{s-1} \exp(-(m+1)x)
    \sum_{j=0}^m \frac{(-\lambda)^j\exp(jx)}{j!}.
\end{eqnarray}
The resulting sum in the integrand is the truncation of the MacLaurin series for $\exp(y)$,
with $y=-\lambda \exp(x)$, i.e.
\begin{eqnarray}
    \label{eq:amore variant generating function c}
    \sum_{j=0}^m \frac{(-\lambda)^j\exp(jx)}{j!}
    = \exp(-\lambda \exp(x))
    -
    \sum_{j=m+1}^\infty \frac{(-\lambda)^j\exp(jx)}{j!}.
\end{eqnarray}
An analysis shows that, for given $\lambda>0$ and as $m \to \infty$, that the
latter sum does not contribute to the main asymptotics of the integral.
Therefore, dropping that sum, and then extending the domain of integration to
$\infty$, we have
\begin{eqnarray}
    \label{eq:j sum amore variant asympt}
    &&\sum_{j=0}^m \frac{(-1)^j}{(j+1)^s \lambda^j (m-j)!} \notag \\
    &\sim& \frac{(-1)^m}{\Gamma(s) \lambda^m} \int_0^\infty x^{s-1} \exp(-\lambda \exp(x)-(m+1)x) dx \notag \\
    &\sim& \frac{(-1)^m \exp(-\lambda)}{\lambda^m (m+1+\lambda)^s},
\end{eqnarray}
with all our estimates above uniform on compact subsets of $\Re{s}>0$.
The last $\sim$ can be
derived by substituting $t=(m+1)x$ into the integral and expanding, for
$t=O(m+1)$,
\begin{equation}
    \lambda \exp(t/(m+1)) = \lambda + \lambda t/(m+1) + O(\lambda t^2/m^2).
\end{equation}

A more careful analysis shows for a given $\lambda$, that to get a good
approximation (i.e. within a specified small relative error $\epsilon$), for a
given $\lambda$, we need $m \gg \lambda$, with the implied constant
depending on $s$ and $\epsilon$.

Formula~\eqref{eq:amore variant thm} becomes interesting as we take $\lambda$ large, because
the terms in that formula are, for $m$ sufficiently large and including the $\exp(-\lambda)$
that appears in front of the sum,
\begin{equation}
    \label{eq:m asympt}
    \sim \frac{(-1)^m \exp(-2\lambda)}{(m+1+\lambda)^s}.
\end{equation}
Thus, while $(-1)^m/(m+1+\lambda)^s$ decays at roughly the same rate as terms in
the Dirichlet series defining $\eta(s)$, our estimate shows that the terms
in~\eqref{eq:amore variant} are eventually, for $m \gg_s \lambda$, exponentially small.
Thus, for $\lambda$ large, there is a transition zone, as we sum over $m$,
in which we quickly achieve roughly $2\lambda/\log(10)$ digits accuracy for $\eta(s)$ before the
sum continues its leisurely convergence to $\eta(s)$. We can also show that the lower term
in~\eqref{eq:m asympt} is, uniformly for given $\lambda>0$ and in compact subsets of $\Re{s}>0$,
$O(1/m^{\Re{s}+2})$ as $m \to \infty$,
thus justifying equation~\eqref{eq:amore variant} in $\Re{s}>0$. Lower terms for this example
will be discussed in a subsequent paper~\cite{R3}.

\subsection{Application to the alternating Hurwitz zeta function}
\label{sec:summation hurwitz}

We define, for $a, \Re{s} >0$,  the alternating Hurwitz zeta function to be
\begin{equation}
    \label{eq:alt hurwitz zeta}
    \zeta^*(s,a) :=
    \sum_{n=0}^{\infty} \frac{(-1)^{n}}{(n+a)^{s}}.
\end{equation}
Setting, in~\eqref{eq:gh summation d}, $g(n)=(-1)^n/(n+a)^{s_0}$ and $h(n)=1/(n+a)^s$, we have
\begin{eqnarray}
    \label{eq:hurwitz G}
    \tilde{G}(z) &=&  \sum_{n=0}^{\infty} \frac{(-1)^{n}}{(n+a)^{s_0}} z^n \notag \\
    &=& \frac{1}{\Gamma(s_0)} \int_0^\infty x^{s_0-1} \sum_{n=0}^\infty (-z)^n \exp(-(n+a)x) dx \notag \\
    &=& \frac{1}{\Gamma(s_0)} \int_0^\infty \frac{x^{s_0-1}\exp(-ax)}{1+z\exp(-x)} dx,
\end{eqnarray}
and so
\begin{eqnarray}
    \label{eq:alt hurwitz zeta b}
    &&\sum_{n=0}^{\infty} \frac{(-1)^{n}}{(n+a)^{s_0}} \frac{1}{(n+a)^s} z^n =\notag \\
    &&\frac{1}{\Gamma(s_0)}
    \sum_{m=0}^\infty z^m
    \left(
        \int_0^\infty \frac{x^{s_0-1}\exp(-x(m+a))}{(1+z\exp(-x))^{m+1}} dx
    \right)
    \sum_{j=0}^m \frac{(-1)^j {m \choose j}}{(j+a)^s}. \notag \\
\end{eqnarray}
We remark that, here, the function $\tilde{G}(z)$ is equal to $\Phi(-z,s_0,a)$,
where $\Phi$ is the Lerch zeta function.

A similar analysis of~\eqref{eq:hurwitz G} to that of~\eqref{eq:integral eta} can
be carried out and shows that~\eqref{eq:alt hurwitz zeta b} holds for $s \in \C$,
$\Re{s_0}>0$ and $z \notin (-\infty,-1]$. The asymptotic formula~\eqref{eq:j sum asympt}
provides a uniform bound for the sum over $j$.

Thus, substituting $z=1$ gives
\begin{eqnarray}
    \label{eq:alt hurwitz zeta c}
    &&\zeta^*(s+s_0,a) = \notag \\
    &&\frac{1}{\Gamma(s_0)}
    \sum_{m=0}^\infty
    \left(
        \int_0^\infty \frac{x^{s_0-1}\exp(-x(m+a))}{(1+\exp(-x))^{m+1}} dx
    \right)
    \sum_{j=0}^m \frac{(-1)^j {m \choose j}}{(j+a)^s}. \notag \\
\end{eqnarray}
We can extend the range of validity to
$\Re{s_0}>-M$ by repeatedly integrating by parts $M$ times. For example,
for $\Re{s_0}>-1$ and $s \in \C$:
\begin{eqnarray}
    \label{eq:alt hurwitz zeta d}
    &&\zeta^*(s+s_0,a) = \notag \\
    &&\frac{1}{\Gamma(s_0+1)}
    \sum_{m=0}^\infty 
    \left(
        \int_0^\infty \frac{x^{s_0}\exp(-x(m+a))(m+a+(a-1)\exp(-x))}{(1+\exp(-x))^{m+1}} dx
    \right) \notag \\
    && \qquad \qquad \qquad \times \sum_{j=0}^m \frac{(-1)^j {m \choose j}}{(j+a)^s}.
\end{eqnarray}
Substituting $s_0=0$, the above integral can easily be evaluated to equal
$1/2^{m+1}$, thus specializing to another formula of Hasse.

We also have, as in Section~\ref{sec:alt zeta}, for all $s_0, s \in \C$,
\begin{equation}
    \label{eq:alt hurwitz interpolate}
    \zeta^*(s+s_0,a) =
    \sum_{m=0}^\infty \frac{(-1)^m}{m!}
    \sum_{l=1}^{m+1} s(m+1,l) \zeta^*(s_0+1-l,a)
    \sum_{j=0}^m \frac{(-1)^j {m \choose j}}{(j+a)^s}.
\end{equation}

We give yet another formula for the alternating Hurwitz zeta function,
when $a=N+1$ is a positive integer.
Taking, in our summation formula~\eqref{eq:gh summation}, $g(n)=0$ if $n\leq N$, and $(-1)^{n-N-1}/n^{s_0}$
if $n \geq N+1$, and $h(n)=1/n^s$, we have:
\begin{eqnarray}
    \label{eq:G alt hurwitz}
    G(z) &=&  \sum_1^\infty g(n) z^{n-1} = (-1)^N \sum_{N+1}^\infty \frac{(-z)^{n-1}}{n^{s_0}} \notag \\
    &=& \frac{z^N}{\Gamma(s_0)} \int_0^\infty x^{s_0-1} \frac{\exp(-(N+1)x)}{1+z\exp(-x)} dx.
\end{eqnarray}
This gives the analytic continuation of $G(z)$ to the point $z=1$ where $G$ has radius of convergence 2.
Therefore, substituting $z=1$ and proceeding as before
\begin{eqnarray}
    \label{eq:alt hurwitz interpolate b}
    &&\zeta^*(s+s_0,N+1) = \notag \\
    &&\sum_{m=0}^\infty \frac{(-1)^m}{m!}
    \sum_{l=1}^{m+1} s(m+1,l) \zeta^*(s_0+1-l,N+1)
    \sum_{j=0}^m \frac{(-1)^j {m \choose j}}{(j+1)^s}. \notag \\
\end{eqnarray}

\subsection{Dirichlet $L$-functions again}
\label{sec:summation dirichlet}

Next we give formulas that express $L(s,\chi)$, for {\it any} non-trivial
character $\chi$ for the modulus $q$ (i.e. without restriction on $q$), in
terms of the Riemann zeta function, or, more precisely, in terms of the
alternating zeta function.

We let $g(n)=(-1)^{n-1}/n^{s_0}$ and $h(n)=(-1)^{n-1} \chi(n)/n^s$. Therefore, $G(z)$ here is identical
to the function in~\eqref{eq:integral eta} and the whole discussion of~\eqref{eq:alt zeta}
regarding the function $G(z)$ carries through. We also require a bound
for the sum over $j$, which we will show beats the
trivial bound exponentially.

Let
\begin{equation}
    \label{eq:C_q defn}
    C_q:=|1+e(1/q)|.
\end{equation}
Assume that $\Re{s}>-M$ where $M$ is a non-negative integer.
We will prove, for $q$ fixed and as $m \to \infty$,
\begin{eqnarray}
    &&\sum_{j=0}^m \frac{\chi(j+1)}{(j+1)^s} {m \choose j} \ll m^M C_q^m,
\end{eqnarray}
uniformly for $s$ on compact subsets of $\Re{s}>-M$ (with the implied constant in the
$\ll$ also depending on $q$ and $M$). But $C_q<2$, hence, dividing by $2^m$ shows that
$\frac{1}{2^m} \sum_{j=0}^m \frac{\chi(j+1)}{(j+1)^s} {m \choose j}$ decreases
exponentially fast as $m \to \infty$. It is crucial for this bound, and also
for getting a nice rate of decay on $G^{(m)}(z)/m!$, that we include here the
factor $(-1)^{n-1}$ in $g(n)$ and $h(n)$. 

Thus,
\begin{equation}
    L(s+s_0,\chi) =
    \sum_{m=0}^\infty (m+1)
    \left(
        \int_0^\infty \frac{x^{s_0}\exp(-x(m+1))}{(1+\exp(-x))^{m+2}} dx
    \right)
    \sum_{j=0}^m \frac{\chi(j+1) {m \choose j}}{(j+1)^s},
    \label{eq:L method 2}
\end{equation}
and
\begin{equation}
    L(s+s_0,\chi) =
    \sum_{m=0}^\infty \frac{(-1)^m}{m!}
    \sum_{l=1}^{m+1} s(m+1,l) \eta(s_0+1-l)
    \sum_{j=0}^m \frac{\chi(j+1) {m \choose j}}{(j+1)^s}.
    \label{eq:L method 2 b}
\end{equation}
Equation~\eqref{eq:L method 2} converges for $\Re{s_0}>-1$, and $s \in
\C$, while~\eqref{eq:L method 2 b} converges for all $s_0,s \in \C$, uniformly
for $s_0,s$ on compact subsets of $\C$, thus also giving the analytic
continuation of $L(s,\chi)$.

\subsection{A bound on $\sum_0^m \chi(j+1) {m \choose j}/(j+1)^s$}

In this section we use the notation
\begin{equation}
    \label{eq:mu}
    e(t) := \exp(2\pi i t).
\end{equation}

We have the following identity
\begin{eqnarray}
    \label{eq:chi sum identity}
    &&\sum_0^m \frac{\chi(j+1)}{(j+1)^s} {m \choose j} \notag \\
    &=&\frac{1}{q}
    \sum_{l=0}^{q-1}
    \sum_{j=0}^m \frac{e(-jl/q)}{(j+1)^s} {m \choose j}
    \sum_{a=1}^{q} \chi(a) e((a-1)l/q),
\end{eqnarray}
which is easily verified by noting that, for given $a$, that only the terms $j$ with
$j+1=a \mod q$ survive the sum over $l$. Each $a$, above, thus contributes
\begin{equation}
    \label{eq:diagonal contribution}
    \sum_{0 \leq j \leq m \atop j+1=a\mod m} \frac{\chi(j+1)}{(j+1)^s} {m \choose j}.
\end{equation}
Summing over all $1\leq a \leq q$ thus gives the lhs of~\eqref{eq:chi sum identity}.

We first examine the $s=0$ case because we can obtain very precise information in that case.
Identity~\eqref{eq:chi sum identity} reads
\begin{eqnarray}
    \label{eq:chi sum s=0}
    &&\sum_{j=0}^m \chi(j+1) {m \choose j} \notag \\
    &=&\frac{1}{q}
    \sum_{l=0}^{q-1}
    (1+e(-l/q))^m
    \sum_{a=1}^{q} \chi(a) e((a-1)l/q),
\end{eqnarray}
Notice, crucially, the sum over $a$ on the rhs vanishes when $l=0$. It is for this reason that
we are able to beat the trivial bound.

Furthermore, $|1+e(-l/q)| < 2$ for
$1 \leq l \leq q-1$, with the largest values being when $l=1,q-1$. Also note
that, above, the sum over $a$ can be expressed in terms of the Gauss sum. When $l=1$
the sum over $a$ equals $e(-1/q) \tau(\chi)$, and when $l=q-1$ the sum equals
$e(1/q) \chi(-1) \tau(\chi)$.

Therefore, for fixed $q$,
\begin{eqnarray}
    &&\sum_{j=0}^m \chi(j+1) {m \choose j} \sim \notag \\
    &&\frac{\tau(\chi)}{q} \left(e(-1/q)(1+e(-1/q))^m+ \chi(-1)e(1/q) (1+e(1/q))^m\right)\notag \\
    &&\leq \frac{2}{q^{1/2}} C_q^m,
    \label{eq:C_q}
\end{eqnarray}
where $C_q=|1+e(1/q)|<2$.
Note that we have also used $|\tau(\chi)|\leq
q^{1/2}$, with equality when $\chi$ is primitive. Dividing by $2^m$, the rhs
decays exponentially in $m$. Therefore, because $G(z)$ here is given
by~\eqref{eq:integral eta} and has Taylor series with radius of convergence 2
at $z=1$, we can substitute $z=1$ in our summation formula thus yielding
equations~\eqref{eq:L method 2} and~\eqref{eq:L method 2 b}, when $s=0$.

We can be a bit more precise about the asymptotics of the middle expression
of~\eqref{eq:C_q}, i.e. of the $l=1,q-1$ terms in~\eqref{eq:chi sum s=0}.
Expand
\begin{equation}
    \label{eq:exp series}
    1+e(t) = 2(1+\pi it -\pi^2 t^2 -\frac{2i}{3}\pi^3t^3+O(t^4)).
\end{equation}
The bracketed term on the rhs of~\eqref{eq:exp series} can be written as
\begin{equation}
    \label{eq:exp series b}
    \exp(\pi it -\pi^2 t^2/2 +O(t^4)).
\end{equation}
Therefore,
\begin{eqnarray}
    &&\frac{\tau(\chi)}{q} \left(e(-1/q)(1+e(-1/q))^m+ \chi(-1)e(1/q) (1+e(1/q))^m\right)\notag \\
    &&= \frac{\tau(\chi)}{q} 2^m e^{-\frac{m\pi^2}{2q^2}(1+O(1/q^2))}
    ( e^{-\frac{\pi i}{q}(m+2)}+\chi(-1)e^{\frac{\pi i}{q}(m+2)} ) \notag \\
    &&= \frac{\tau(\chi)}{q} 2^{m+1} e^{-\frac{m\pi^2}{2q^2}(1+O(1/q^2))}
    \times
    \begin{cases}
         \cos({\frac{\pi}{q}(m+2)}), &\text{if $\chi(-1)=1$,} \\
        -i\sin({\frac{\pi}{q}(m+2)}), &\text{if $\chi(-1)=-1$.} \notag \\
    \end{cases} \\
    \label{eq:1,q-1 asympt}
\end{eqnarray}
Note that in the case that $\chi$ is a real primitive character, the above simplifies further to
\begin{eqnarray}
    \frac{2^{m+1}}{q^{1/2}} e^{-\frac{m\pi^2}{2q^2}(1+O(1/q^2))}
    \times
    \begin{cases}
         \cos({\frac{\pi}{q}(m+2)}), &\text{if $\chi(-1)=1$,} \\
         \sin({\frac{\pi}{q}(m+2)}), &\text{if $\chi(-1)=-1$.} \notag \\
    \end{cases} \\
\end{eqnarray}

We can similarly bound the contribution from the terms $2\leq l \leq q-2$ in~\eqref{eq:chi sum s=0}.
For simplicity assume that $\chi$ is primitive so that the inner sum in~\eqref{eq:chi sum s=0}
can be written in terms of $\tau(\chi)$ and is, in absolute value, equal to $q^{1/2}$ or 0,
depending whether $\gcd(l,q)=1$ or not.
Furthermore, amongst $2\leq l \leq q-2$, the maximum $|1+e(l/q)|^m$ occurs when $l=2$ or $q-2$.
Thus the contribution to~\eqref{eq:chi sum s=0} from  $2\leq l \leq q-2$ is bounded by
\begin{equation}
    2^m q^{1/2}
    e^{-\frac{2m\pi^2}{q^2}(1+O(1/q^2))}.
\end{equation}
Notice that the exponent here is four times larger than in~\eqref{eq:1,q-1 asympt}.

We have thus shown, for primitive $\chi$, that
\begin{eqnarray}
    &&\sum_0^m \chi(j+1) {m \choose j} \notag \\
    &&= \frac{\tau(\chi)}{q} 2^{m+1} e^{-\frac{m\pi^2}{2q^2}(1+O(1/q^2))}
    \times
    \begin{cases}
         \cos({\frac{\pi}{q}(m+2)}),  &\text{if $\chi(-1)=1$,} \\
        -i\sin({\frac{\pi}{q}(m+2)}), &\text{if $\chi(-1)=-1$.} \notag \\
    \end{cases} \\
    &&+R(\chi,m),
    \label{eq:chi j sum asympt}
\end{eqnarray}
where the remainder $R(\chi,m)$ is bounded, in absolute value, by
\begin{equation}
    2^m q^{1/2} e^{-\frac{2m\pi^2}{q^2}(1+O(1/q^2))}.
\end{equation}
With respect to uniform asymptotics allowing both $q$ and $m$ to grow,
equation~\eqref{eq:chi j sum asympt} provides
an asymptotic formula for the sum on the left as $m/q^2
\to \infty$ and $m/q^4 \to 0$,
\begin{eqnarray}
    &&\sum_{j=0}^m \chi(j+1) {m \choose j} \notag \\
    &&\sim \frac{\tau(\chi)}{q} 2^{m+1} e^{-\frac{m\pi^2}{2q^2}}
    \times
    \begin{cases}
         \cos({\frac{\pi}{q}(m+2)}), &\text{if $\chi(-1)=1$,} \\
        -i\sin({\frac{\pi}{q}(m+2)}), &\text{if $\chi(-1)=-1$,} \notag \\
    \end{cases} \\
\end{eqnarray}
while the second line of~\eqref{eq:C_q} provides an asymptotic formula valid in
the larger region $m/q^2 \to \infty$. Note that we do no need to assume that
$\chi$ is primitive in the above asymptotic formula- that assumption was made
to simplify the discussion regarding $R(\chi,m)$, but a similar analysis
applies for a general non-trivial character.

We interpret the above asymptotic formula to be an equality in the event that
the rhs vanishes. There are two cases to consider. If $\chi(-1)=-1$, then the
rhs vanishes iff $m = -2 \mod q$. On the other hand, the terms $j$ and $m-j$ on
the lhs cancel each other out, assuming that $j\neq m-j$, because the binomial
coefficients match, and, furthermore, when $m = -2 \mod q$,
$\chi(m-j+1)=\chi(-2-j+1) =\chi(-1)\chi(j+1)=-\chi(j+1)$. Special care is
needed when $m$ is even and $j=m/2$ because then there is just one binomial
coefficient, and not a pair. But, if $j=m/2$ then $\chi(j+1) = \chi((m+2)/2)$.
However, $m+2$ is assumed, here, to be a multiple of $q$, hence $\gcd((m+2)/2,q)>1$ (because $q>2$
when $\chi$ is non-trivial), and so $\chi((m+2)/2)=0$. Thus, $\chi(j+1)=0$, and
the middle binomial coefficient does not contribute to the sum.

The case $\chi(-1)=1$ is a bit more complicated, but the idea is the same,
namely to show that the terms $j$ and $m-j$ cancel. We have
$\cos({\frac{\pi}{q}(m+2)})=0$ iff $m+2 =rq/2$ where $r$ is an odd integer,
i.e. $2m+4=rq$. Thus $q$ is even because $r$ is odd, and, in fact, divisible by
at least $4$ because there are no even primitive characters for the modulus $2$.
Hence
$m=-2 \mod q/2$. Because $r$ is odd, we can write it as $r=2r_0+1$. Now,
$\chi(m-j+1)= \chi((2r_0+1)q/2-j-1)=\chi(q/2-j-1)=\chi(j+1-q/2)$, the last step
because $\chi(-1)=1$.

Next, factor $q$ into distinct prime powers,
$q=p_1^{\beta_1} p_2^{\beta_2} \ldots p_k^{\beta_k}$, with $p_1=2$, and let
$\chi(n)=\chi_{p_1^{\beta_1}}(n) \ldots\chi_{p_k^{\beta_k}}(n)$. Now, for
$1<l\leq k$, $\chi_{p_l^{\beta_l}}(j+1-q/2)= \chi_{p_l^{\beta_l}}(j+1)$, because
$q/2=0 \mod p_l^{\beta_j}$ for $l>1$. Thus, to prove here that $\chi(m-j+1) =
-\chi(j+1)$, i.e. that the terms $j$ and $m-j$ cancel, we need to show that
\begin{equation}
    \label{eq:chi power 2 goal}
    \chi_{2^{\beta_1}}(j+1-q/2) = -\chi_{2^{\beta_1}}(j+1)
\end{equation}
for all $j$ and any
primitive character mod $2^{\beta_1}$.

Note that we can assume that $j+1$ is odd
because $q/2$ is even. Thus, if $j+1$ is even then
both $\chi_{2^{\beta_1}}(j+1-q/2)$ and $\chi_{2^{\beta_1}}(j+1)$ are equal to 0.

Write $q=2^{\beta_1}(2q_0+1)$, i.e. write the odd part of $q$ as $2q_0+1$.
Thus, $\chi_{2^{\beta_1}}(j+1-q/2) = \chi_{2^{\beta_1}}(j+1-2^{\beta_1-1})$.

Now, if $\beta_1=2$, then $\chi_4(j+1)=-\chi_4(j-1)$ and we are done.
Next, if $\beta_1\geq3$ then any character for the
modulus $2^{\beta_1}\geq 8$ can be written, for odd $n$, as
\begin{equation}
    \label{eq:chi power 2}
    \chi_{2^{\beta_1}}(n) = e(\frac{\mu\nu}{2}+\frac{\mu'\nu'}{2^{\beta_1-2}})
\end{equation}
where $0 \leq \mu<2$, $0\leq \mu' < 2^{\beta_1-2}$, and
\begin{equation}
    \label{eq:n}
    n \equiv (-1)^\nu 5^{\nu'} \mod 2^{\beta_1},
\end{equation}
with $0 \leq \nu<2$, $0\leq \nu' < 2^{\beta_1-2}$
See Chapters 4-5 of~\cite{D}.
Furthermore, $\chi_{2^{\beta_1}}$ is primitive implies that $\mu'$ is odd.

We claim that
\begin{equation}
    \label{eq:claim 1}
    n-2^{\beta_1-1} \equiv n(2^{\beta_1-1}+1) \mod 2^{\beta_1}
\end{equation}
which holds because $n$ is assumed to be odd, and
\begin{equation}
    \label{eq:claim 2}
    5^{2^{\beta_1-3}} \equiv 2^{\beta_1-1}+ 1 \mod 2^{\beta_1},
\end{equation}
which can be proven inductively by repeatedly squaring.

Therefore
\begin{equation}
    \label{eq:almost there}
    \chi_{2^{\beta_1}}(n-2^{\beta_1-1}) = \chi(n) \chi(2^{\beta_1-1}+1) = \chi(n) \chi(5^{2^{\beta_1-3}}).
\end{equation}
But
\begin{equation}
    \label{eq:chi almost there}
    \chi_{2^{\beta_1}}(5^{2^{\beta_1-3}}) = e(\frac{\mu'2^{\beta_1-3}}{2^{\beta_1-2}}) = \exp(\pi i \mu') = -1,
\end{equation}
the last step because $\mu'$ is odd when $\chi_{2^{\beta_1}}$ is primitive.

We have therefore shown that if: $\chi$ is primitive, $\chi(-1)=1$, and
$\cos({\frac{\pi}{q}(m+2)})$$=0$, then $\chi(j+1) = -\chi(m-j+1)$ for all $j$.
Finally, we need to address, as before, the case where $m$ is even and $j=m/2$,
i.e. when the middle binomial coefficient is not part of pair. But, in that
case, $\chi(j+1)= \chi((m+2)/2) = \chi(rq/2)$ which equals 0 because $q$ is
divisible by $4$ and hence $q/2$ is even.

\subsubsection{Bounding the sum in general}

Assume for now that $\Re{s}>0$.

We consider the sum over $j$ in~\eqref{eq:chi sum identity}, which we can write as
\begin{eqnarray}
    \label{eq:chi j sum}
    &&\sum_{j=0}^m \frac{e(-jl/q)}{(j+1)^s} {m \choose j} \notag \\
    &=& \frac{1}{\Gamma(s)}
    \int_0^\infty x^{s-1} (1+\exp(-x) e(-l/q))^m \exp(-x) dx.
\end{eqnarray}
Similar to the $s=0$ case, one can show that the $l=1,q-1$ terms dominate (the $l=0$
term in~\eqref{eq:chi sum identity} vanishes), and we have, for fixed $q$ and as $m \to \infty$,
\begin{eqnarray}
    &&\sum_{j=0}^m \frac{\chi(j+1)}{(j+1)^s} {m \choose j} \sim \notag \\
    &&\frac{\tau(\chi)}{q\Gamma(s)} \bigg(
    e(-1/q)\int_0^\infty x^{s-1} (1+\exp(-x) e(-1/q))^m \exp(-x) dx \notag \\
    &+& \chi(-1)e(1/q) \int_0^\infty x^{s-1} (1+\exp(-x) e(1/q))^m \exp(-x) dx \bigg). \notag \\
    \label{eq:chi sum s}
\end{eqnarray}
Taking absolute value and using $|1+\exp(-x) e(-1/q)| \leq C_q$
gives~\eqref{eq:chi sum bound any s}

Repeated integration by parts of~\eqref{eq:chi j sum}, say $M$ times, allows us
to obtain the same sort of bound for $\Re{s}>-M$, at a cost of, on integrating
by parts, introducing extra powers of $m$ into the integrand, thus giving
bound~\eqref{eq:chi sum bound any s}.

This rate of decay, beating
$2^m$ by an exponentially small amount, justifies formulas~\eqref{eq:L method 2}
and~\eqref{eq:L method 2 b}.


\section{An estimate for $\sum_{j=0}^m \frac{(-1)^j}{(j+a)^s} {m \choose j}$}

For $a>0$, we derive an asymptotic formula for
\begin{eqnarray}
    \label{eq:j sum a}
    \sum_{j=0}^m \frac{(-1)^j}{(j+a)^s} {m \choose j},
\end{eqnarray}
as $m \to \infty$.
The above equals, for $m+\Re{s}>0$,
\begin{eqnarray}
    \label{eq:derive asympt j sum}
    &&\frac{1}{\Gamma(s)}
    \int_0^\infty x^{s-1} \exp(-ax) (1-\exp(-x))^m dx = \notag \\
    &&\frac{1}{m^a \Gamma(s)}
    \int_0^m (-\log(v/m))^{s-1} v^{a-1} (1-v/m)^{m} dv.
\end{eqnarray}
Here we have substituted $v=m\exp(-x)$.
Now, we can write
\begin{equation}
    \log(m/v)^{s-1} = \log{m}^{s-1} ( (1-\log{v}/\log{m})^{s-1}-1 + 1),
\end{equation}
so that~\eqref{eq:derive asympt j sum} equals
\begin{eqnarray}
    \label{eq:derive asympt j sum b}
    &&\frac{\log(m)^{s-1}}{m^a \Gamma(s)}
    \int_0^m v^{a-1} (1-v/m)^{m} dv + \notag \\
    &&\frac{\log(m)^{s-1}}{m^a \Gamma(s)}
    \int_0^m ( (1-\log{v}/\log{m})^{s-1}-1) v^{a-1} (1-v/m)^{m} dv. \notag \\
\end{eqnarray}
It is a simple exercise to show that the first integral above is
asymptotically equal to $\Gamma(a)$ as $m \to \infty$ by breaking up the
integral into, say, $\int_0^{m^{1/4}}+\int_{m^{1/4}}^m$, and using:
\begin{eqnarray}
    &&(1-v/m)^m = \exp(m\log(1-v/m)) \notag \\
    &&=\exp(-v-\frac{v^2}{2m}-\frac{v^3}{3m}-\ldots), \qquad |v|<m \notag \\
    &&=\exp(-v)(1+O(1/m^{1/2})), \qquad |v| < m^{1/4}.
\end{eqnarray}
We can also show that the second integral in~\eqref{eq:derive asympt j sum b}
does not contribute to the main term.
First, undo our change of variable and express the integral as
\begin{equation}
    \label{eq:remainder int}
    \int_0^\infty ( (x/\log{m})^{s-1}-1) (m\exp(-x))^{a} (1-\exp(-x))^{m} dx. \notag \\
\end{equation}
We need to show that this tends to 0 as $m \to \infty$, and uniformly for $s$ in compact
subsets of $\C$. To do so,
break the above into three integrals:
\begin{eqnarray}
    \label{eq:three integrals}
    \int_0^\infty &=& 
    \int_0^{\log(m) - \log(m)^{1/2}}
    \notag \\
    &+&\int_{\log(m) - \log(m)^{1/2}}^{\log(m) + \log(m)^{1/2}}
    +\int_{\log(m) + \log(m)^{1/2}}^\infty.
\end{eqnarray}
In the middle integral, we have $(x/\log{m})^{s-1}-1=O(\log(m)^{-1/2})$, with
the implied constant uniform for $s$ in compact subsets of $\C$.
Hence, pulling this estimate out, and using, as before,
\begin{equation}
    \label{eq:remainder int b}
    \int_0^\infty (m\exp(-x))^{a} (1-\exp(-x))^{m} dx \sim \Gamma(a)
\end{equation}
we have
\begin{eqnarray}
    \label{eq:middle integral}
    &&\int_{\log(m) - \log(m)^{1/2}}^{\log(m) + \log(m)^{1/2}} 
    ( (x/\log{m})^{s-1}-1) (m\exp(-x))^{a} (1-\exp(-x))^{m} dx \notag \\
    &=&  O(\log(m)^{-1/2}).
\end{eqnarray}

Next we bound the first integral, using
\begin{eqnarray}
    &&(1-\exp(-x))^m =\notag \\
    &&\exp(-m(\exp(-x)+\exp(-2x)/2+\exp(-3x)/3+\ldots))
    \notag \\
    &&\leq \exp(-m\exp(-x)).
    \label{eq:second integral bound a}
\end{eqnarray}
Note, further, that $(1-\exp(-x))^m$ increases monotonically for $x\geq 0$, which, combined
with the above estimate gives
\begin{equation}
    (1-\exp(-x))^m \leq \exp(-\exp(\log(m)^{1/2})),
    \label{eq:second integral bound b}
\end{equation}
when $0 \leq x \leq \log(m)-\log(m)^{1/2}$. This shows that $(1-\exp(-x))^m$
decreases much faster than any power of $m$, in particular $m^{a}$.
Furthermore, if $\Re{s}\geq 1$, then $(x/\log{m})^{s-1}-1$ is bounded over
the interval $[0,\log(m) - \log(m)^{1/2}]$. If $\Re{s}<1$, then let $\Re{s}>-M$, where $M$ is
a non-negative integer. For $m>M$, we can steal $M+1$ powers of $(1-\exp(-x))=x+O(x^2)$,
for $|x|\leq 1$, from $(1-\exp(-x))^m$, with little effect on that factor, so that
$(x/\log{m})^{s-1}-1)(1-\exp(-x))^{M+1}$ is, on the interval $[0,1]$ of size
$O(\log(m)^{-\sigma+1} +1)$, where $\sigma = \Re{s}$. Furthermore, the latter $O$ estimate also holds
for $(x/\log{m})^{s-1}-1)$
once $x>1$, if $\Re{s}<1$. The integrand of the second integral is thus handily dominated by the
last factor of size $O(\exp(-\exp(\log(m)^{1/2})))$ and the second integral tends to $0$,
uniformly for $s$ in compact subsets of $\C$.

Finally, to bound the third integral, consider the factor $(x/\log{m})^{s-1}-1$.
If $\Re{s}\leq 1$ then this is bounded for $x\geq \log(m) + \log(m)^{1/2}$.
And if $\Re{s}>1$ then this factor is $O(x^{\sigma-1}/\log(m)^{\sigma-1})$.
Using the trivial bound $(1-\exp(-x))^m < 1$, we thus have an upper bound for
the third integral, when $\Re{s} \leq 1$, of
\begin{eqnarray}
     \ll m^a \int_{\log(m)+\log(m)^{1/2}}^\infty \exp(-ax) dx
     = \frac{\exp(-a \log(m)^{1/2})}{a}
    \label{eq:bound third integral}
\end{eqnarray}
and, when $\Re{s}>1$, of
\begin{eqnarray}
    \ll \frac{m^a}{\log(m)^{\sigma-1}} \int_{\log(m)+\log(m)^{1/2}}^\infty x^{\sigma-1} \exp(-ax) dx.
    \label{eq:bound third integral b}
\end{eqnarray}
Substituting $u=x-\log(m)-\log(m)^{1/2}$, then pulling out
$(\log(m)+\log(m)^{1/2})^{\sigma-1}$ from the integral, and finally using
$1+y<e^y$, with $y=u/(\log(m)+\log(m)^{1/2})$, the above is bounded by
\begin{eqnarray}
    &&\exp(-a\log(m)^{1/2}) \frac{(\log(m)+\log(m)^{1/2})^{\sigma-1}}{\log(m)^{\sigma-1}}
    \notag \\
    &&\times \int_{0}^\infty \exp(u(-a+\frac{(\sigma-1)}{\log(m)+\log(m)^{1/2}})) du \notag \\
    &&\ll  \exp(-a\log(m)^{1/2}),
    \label{eq:bound third integral c}
\end{eqnarray}
with the implied constant, for given $a$,  uniform for $s$ in compact subsets of $\C$.

We have thus shown, for given $a$, that the second integral
in~\eqref{eq:derive asympt j sum b} tends
uniformly, for $s$ in compact subsets of $\C$, to 0, hence proving
the asymptotic formula~\eqref{eq:j sum asympt}.

Note that when $s \in \Z$ and $s\leq 0$, then the first line in~\eqref{eq:derive asympt j sum},
valid for $m>-s$, vanishes because $1/\Gamma(s)$ is then equal to $0$. Therefore, the asymptotic
formula~\eqref{eq:j sum asympt} should be interpreted as an equality for such $s$ and $m$.

\bibliographystyle{amsplain}


\end{document}